\documentclass[3p,11pt]{elsarticle}
\usepackage{amsmath}
\usepackage{graphicx}
\usepackage{wrapfig}
\usepackage{amsmath}
\usepackage{times}
\usepackage{amssymb}
\usepackage{array}
\usepackage{latexsym}
\usepackage{amsthm}
\usepackage{amsfonts}
\usepackage{nomencl}

\makenomenclature
\usepackage{hyperref}
\usepackage{array,tabularx}
\usepackage{scalerel,stackengine}
\usepackage{graphics}
\usepackage[english]{babel}
\allowdisplaybreaks

\stackMath
\numberwithin{equation}{section}
\date{}
\newtheorem{t1}{Theorem}[section]
\newtheorem{p1}{Proposition}[section]
\newtheorem{l1}{Lemma}[section]

\newtheorem{df}{Definition}[section]

\usepackage[utf8]{inputenc}
\begin{document}
\begin{frontmatter} 
\title{Existence of solutions to gas expansion problem through a sharp corner for 2-D Euler equations with general equation of state}

\author{Rahul Barthwal}
\author{T. Raja Sekhar}
\address{Department of Mathematics, Indian Institute of Technology Kharagpur, Kharagpur,  India} 

%% or include affiliations in footnotes: 
%\author[mymainaddress,mysecondaryaddress]{Elsevier Inc} 
%\ead[url]{www.elsevier.com} 

%\author[mysecondaryaddress]{Global Customer Service\corref{mycorrespondingauthor}} 
\cortext[mycorrespondingauthor]{Corresponding author} 
\ead{trajasekhar@maths.iitkgp.ac.in} 

%\address[mymainaddress]{1600 John F Kennedy Boulevard, Philadelphia} 
%\address[mysecondaryaddress]{360 Park Avenue South, New York} 

\begin{abstract} 
In this article, we study the gas expansion problem by turning a sharp corner into vacuum for the two-dimensional pseudo-steady compressible Euler equations with a convex equation of state. This problem can be considered as interaction of a centered simple wave with a planar rarefaction wave. In order to obtain the global existence of solution up to vacuum boundary of the corresponding two-dimensional Riemann problem, we consider several Goursat type boundary value problems for 2-D self-similar Euler equations and use the ideas of characteristic decomposition and bootstrap method. Further, we formulate two-dimensional modified shallow water equations newly and solve a dam-break type problem for them as an application of this work. Moreover, we also recover the results from the available literature for certain equation of states which provide a check that the results obtained in this article are actually correct. 
\end{abstract} 
\end{frontmatter} 
%\linenumbers
%%%%%%%%%% Insert the texts which can accomdate on fi

\section{Introduction}
\begin{figure}
    \centering
    \includegraphics[width=3 in]{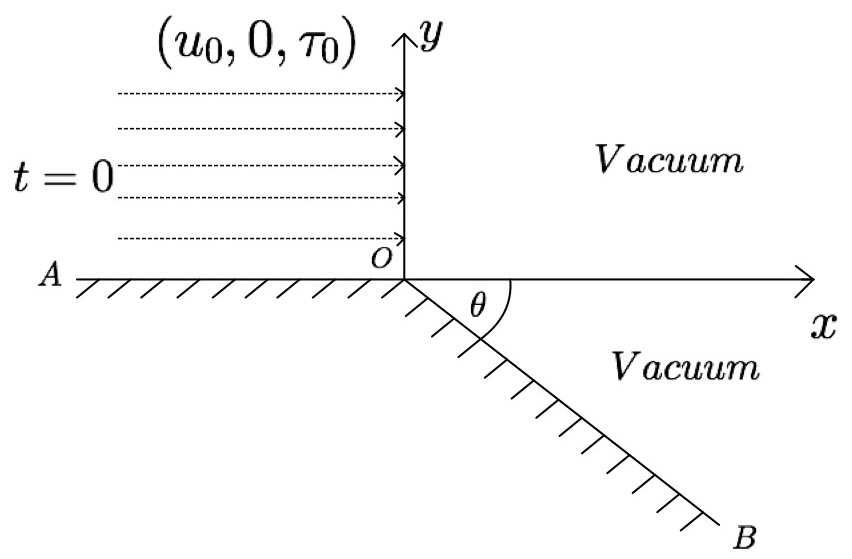}
    \caption{Initial and boundary conditions}
    \label{fig: 1}
\end{figure}

The mathematical theory of compressible flows in two dimensions developed rapidly in recent years. Supersonic flow around a sharp corner is one of the most important and well-studied elementary flows in the study of compressible flows. Courant and Friedrichs \cite{courant0} noted that the supersonic flow around a bend or sharp corner is affected by simple waves (compression or expansion wave). Consider an infinitely long wedge $OB$ with a horizontal wall $AO$ which is straight up to a sharp corner $O$. Initially, a supersonic flow arrives with a constant velocity $(u_0,0)$ and density $\rho_0$ along the ground wall and then suddenly expands to vacuum in the other region of the corner; see Figure \ref{fig: 1}.  The turn of the gas from the peak $O$ is affected by a centered simple wave which starts interacting with a planar rarefaction wave. Therefore, the problem of gas expansion can be essentially considered as the interaction of a centered rarefaction wave with a planar rarefaction wave. The gas expansion problem through a sharp corner has been well-studied recently by a group of mathematicians. Sheng and You \cite{sheng2018interaction} studied this problem for isentropic Euler equations with polytropic gas and proved the existence of a global solution in the entire interaction region. Chen et al. \cite{chen2020expansion} extended their ideas to magnetogasdynamics system with polytropic gas and established an existence result for a global solution. Recently, Lai and Sheng \cite{lai2021two} also studied the polytropic gas expansion problem for 2-D Euler equations by turning the gas around a sharp corner into a vacuum for different cases of the wall angles and obtained some beautiful results. An obvious question that may arise in a reader's mind after the success of these works is that can one generalize these results to any arbitrary equation of state. Inspiring by this idea, in this work we are trying to generalize these results for any general convex equation of state. For simplicity we assume that the inclination angle $\theta$ of the wall $OB$ satisfies $\theta\in(-\pi/2, 0)<\alpha_v$, where $\alpha_v$ denotes the $C_+$ characteristic angle when the flow arrives at the vacuum state. 

Let us consider the two-dimensional isentropic Euler equations of gas dynamics \cite{li2011semi}
\begin{equation}\label{eq: 1}
\begin{aligned}
\rho_t+(\rho u)_x+(\rho v)_y=0,\\
(\rho u)_t+(\rho u^2+p)_x+(\rho u v)_y=0,\\
(\rho v)_t+(\rho u v)_x+(\rho v^2+p)_y=0,
\end{aligned}
\end{equation}
where $\rho$ denotes the density, $u$ and $v$ denotes the flow velocity in the $x$ and $y$ direction, respectively, $p=p(\tau)$ denotes the pressure of the gas and $\tau=1/\rho$ is the specific volume. For simplicity of notation and computations, we use $(u, v, \tau)$ as our primitive variables instead of $(u, v, \rho)$.

Cauchy problem for the system \eqref{eq: 1} is a complicated and challenging open problem. The two-dimensional Riemann problem is a particular kind of Cauchy problem which consists of constant initial data along any ray passing through the origin. The expansion problems of a flow into a vacuum are usually special cases of two-dimensional Riemann problems, which are concerned with the interaction of planar rarefaction waves and/or centered rarefaction waves. In recent years, a lot of significant work has been done for the two-dimensional compressible Euler system as well as numerous other related models for a variety of initial and boundary value problems; see viz. \cite{barthwal2022existence, BARTHWAL2023127022, chen2019semi, li1998two, li2011semi, song2009semi, zheng2012systems, barthwal2022simple, hu2018improved, hu2019regularity, barthwalconstruction, barthwal2022two}. In particular for gas expansion problems through a sharp corner or wedge, we refer the reader to \cite{lai2015expansion,  lai2019global, lai2021interactions, li2002two, li2011characteristic, sheng2018two, zafar2016expansion}  and references cited therein. The study of two-dimensional Riemann problems is significant in theoretical and numerical analysis and many engineering applications too; see \cite{glimm2008transonic,kurganov2002solution,lax1998solution, tesdall2007triple}. So the study in this article is of utter importance.
We refer Figure \ref{fig: 1} to impose an initial data on the system \eqref{eq: 1} of the form
\begin{align}\label{eq: 2}
    (u, v, \tau)(x, y, 0)=\begin{cases}
        (u_0, 0, \tau_0), \hspace{1 cm}x<0, y> 0,\\
        \mathrm{vacuum},\hspace{1.25 cm} x>0,y\geq x\tan \theta,
    \end{cases}
\end{align}
where $u_0, \tau_0$ are constants. Clearly, \eqref{eq: 1} with the initial data \eqref{eq: 2} is a 2-D Riemann problem with a boundary. Our main objective in this article is to solve this problem for any arbitrary convex pressure.

Throughout the article we assume that the pressure $p(\tau)$ satisfies the following properties:
\begin{align}
    p'(\tau)<0, p''(\tau)>0 ~~\mathrm{for}~\tau>\tau_0,
\end{align}
which is generally true for most of the cases of the equation of states of physical relevance.

One of the major difficulties in establishing the global existence of solution for the system \eqref{eq: 1} is that the system \eqref{eq: 1} may change its type from hyperbolic to elliptic in the interaction domain and the type of the system is not a priori known. Since different types of partial differential equations involve different solving notions, we need to establish {\it a priori} estimate of the physical variables. We can not use the method of characteristics in the elliptic domain, therefore, to skip the possible occurrence of a bad case of mixed type, we need to use the ideas of characteristic decompositions and invariant regions \cite{li2009interaction} to maintain the hyperbolicity of the system \eqref{eq: 1} in the interaction domain. The other complexity of this article is to handle a general convex equation of state when compared to a polytropic equation of state since many important characteristic functions of density, which are crucial in developing a priori bounds of solutions, may change their behaviour from increasing to decreasing or vice versa in the interaction domain as we will see in Section 4. 

The rest of the article is organized as follows. We reduce system \eqref{eq: 1} in the form of self-similar variables and obtain characteristic decompositions of density and characteristic angles in Section \ref{2} which are helpful for developing \textit{a priori} estimates of physical variables. In Section \ref{3}, we provide expressions for planar rarefaction wave and centered rarefaction wave and establish the boundary data  estimates to prove the existence of a local solution. We construct the invariant regions for characteristic angles and obtain the $C^0$ and $C^1$ norm estimates of the physical variables in the interaction region in Section \ref{4}. Section \ref{5} is devoted to prove the existence of a global  solution by extending the local solution up to the vacuum boundary by solving several Goursat problems locally in each extension step. Further,  we discuss some particular cases of equations of states to discuss some relevant physical models as applications of this work in Section \ref{6}. In particular, we formulate, first time in the literature, two-dimensional modified shallow water equations to solve a dam-break type problem and also recover results of gas expansion problem for certain equations of states such as the polytropic equation of state from the available literature. In Section \ref{7} we finally provide the concluding remarks and future scope of this work.
\section{System in $(\xi, \eta)$ plane}\label{2}
It is easy to see that the system \eqref{eq: 1} and initial data \eqref{eq: 2} are invariant under the transformation $(t, x, y)\longrightarrow (\alpha t, \alpha x, \alpha y)$ for $\alpha>0$. Then we can reduce the Euler system \eqref{eq: 1} in self-similar co-ordinates $\left(\xi=\dfrac{x}{t}, \eta=\dfrac{y}{t}\right)$ as follows 
\begin{equation}\label{eq: 2.1}
\begin{aligned}
     (\rho U)_\xi+ (\rho V)_\eta+2\rho=0,\\
    U U_\xi+ V U_\eta+ \tau p_\xi+ U=0,\\
        U V_\xi+ V U_\eta+ \tau p_\eta+ V=0,
\end{aligned}
\end{equation}
where $U=u-\xi$ and $V=v-\eta$ denotes the components of pseudo-velocity in $(\xi, \eta)$ plane.

Also, the initial data \eqref{eq: 2} now changes into 
\begin{align}\label{eq: 2.2}
    (u, v, \tau)(\xi, \eta)=\begin{cases}
    (u_0, 0, \tau_0), ~~~~~~\xi<0, \eta>0, \\
    \mathrm{vacuum}, ~~~~~~~~\xi>0, \eta\geq \xi \tan \theta, \xi^2+\eta^2\longrightarrow \infty. 
    \end{cases}
\end{align}

Assuming that the flow is irrotational, one can introduce a function $\phi(\xi, \eta)$ such that $U\phi_\xi+V\phi_\eta=q^2$, where $q=\sqrt{U^2+V^2}$. The function $\phi$ is usually referred to as a potential function; see \cite{li2011characteristic}. Further, using the last two equations of system \eqref{eq: 2.1} it is easy to get the pseudo-Bernoulli's law of the form
\begin{equation}\label{eq: 2.3}
\dfrac{q^2}{2}+\displaystyle{\int_{\tau_0}^{\tau}\tau p'(\tau) d\tau}+\phi= 0.
\end{equation}
Under the assumption that the flow is irrotational and for a smooth solution, \eqref{eq: 2.1} can be reduced into a matrix form as follows
\begin{equation}\label{eq: 2.4}
\begin{aligned}
    \begin{bmatrix}
    u \\ v
    \end{bmatrix}_\xi+\begin{bmatrix}
     \dfrac{-2UV}{c^2-U^2} & \dfrac{c^2-V^2}{c^2-U^2}\\ -1 & 0
    \end{bmatrix} \begin{bmatrix}
    u \\ v
    \end{bmatrix}_\eta=0,
\end{aligned}
\end{equation}
where $c(\tau)=\sqrt{-\tau^2 p'(\tau)}$ is the speed of sound.
 
It is straightforward to see that the eigenvalues of the system \eqref{eq: 2.4} are $\lambda_\pm=\dfrac{UV\pm c\sqrt{U^2+V^2-c^2}}{U^2-c^2}$ with corresponding left eigenvectors $l_\pm=(1, \lambda_\mp)$. The expression of these eigenvalues shows that the system \eqref{eq: 2.4} is a mixed type system and changes its behaviour from hyperbolic to elliptic across the sonic boundary and depends on the choice of pseudo-Mach number $M=\dfrac{\sqrt{U^2+V^2}}{c}$. For $M> 1$(supersonic) system \eqref{eq: 2.4} is hyperbolic while for $M<1$(subsonic) it is elliptic. Then we define the two families of wave characteristics as 
\begin{equation}
    \dfrac{d\eta}{d\xi}=\lambda_\pm.
\end{equation}
We call these characteristic curves positive ($C_+$) and negative ($C_-$) characteristics, respectively.

Moreover, we obtain the characteristic equations by multiplying $l\pm$ to the system $\eqref{eq: 2.4}$ as 
\begin{align}\label{eq: 2.6}
\begin{cases}
\bar{\partial}_{+}u+\lambda_- \bar{\partial}_{+}v=0,\\
\bar{\partial}_{-}u+\lambda_+ \bar{\partial}_{-}v=0
\end{cases}
\end{align} 
where $\bar{\partial}_{\pm}=\partial_{\xi}+\lambda_\pm \partial_{\eta}$.

Now we define the concept of characteristic angles as in \cite{li2011characteristic}. The $C_+$ characteristic angle $\alpha$ is defined as the angle between the $C_+$ characteristic direction and $\xi$-axis. In a similar manner, one can define the $C_-$ characteristic angle $\beta$. It is trivial to see that the eigenvalues $\lambda_\pm$ satisfy $\tan \alpha=\lambda_+, \tan \beta=\lambda_-$. Moreover, if we denote pseudo-Mach angle by $\delta$ and pseudo-flow angle by $\sigma$ then $\sigma=\dfrac{\alpha+\beta}{2},~~\delta=\dfrac{\alpha-\beta}{2}$, where $\delta$ is the angle between $C_+(C_-)$ characteristic and pseudo-velocity vector $(U, V)$ and $\sigma$ is the angle between $(U, V)$ and $\xi$-axis; see Figure \ref{fig: 2}. Therefore, we have the relations of the form \cite{li2011characteristic}
\begin{align}\label{eq: 2.7}
\begin{cases}
    u-\xi=c \dfrac{\cos \sigma}{\sin \delta},\\
    v-\eta=c \dfrac{\sin \sigma}{\sin \delta},\\
\sin \delta=\dfrac{1}{M}.
\end{cases}
\end{align}
\begin{figure}
    \centering
    \includegraphics[width=3.5 in]{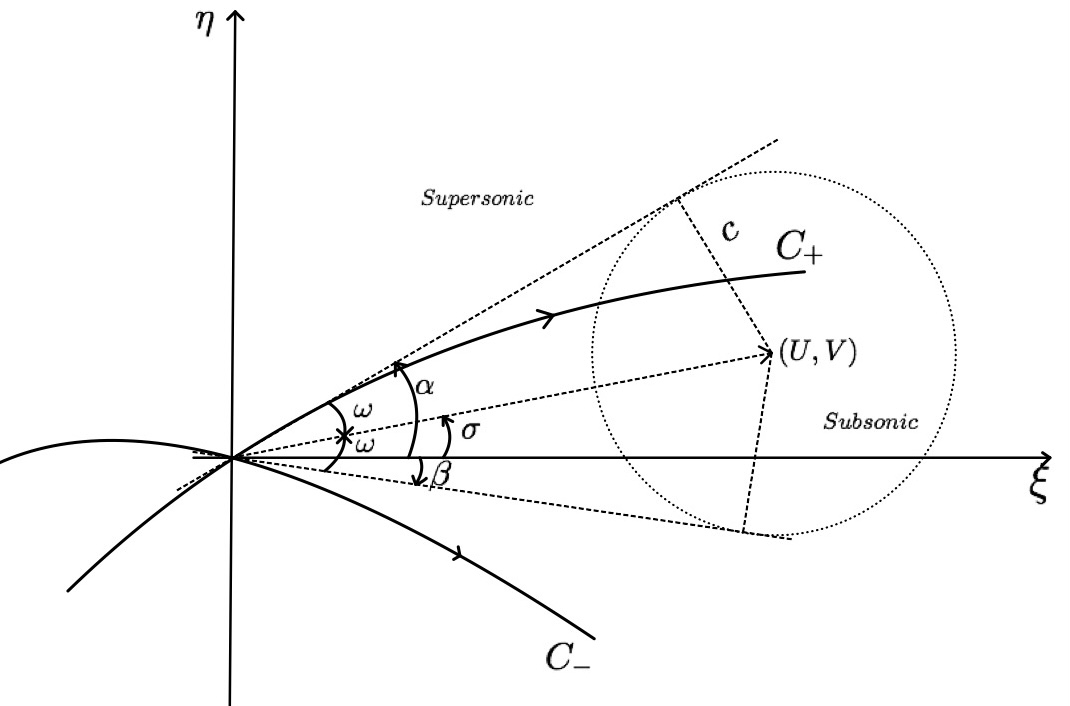}
    \caption{Characteristic angles and pseudo-flow directions}
    \label{fig: 2}
\end{figure} 
\subsection{First order characteristic decompositions}
In this subsection, we derive first-order characteristic decompositions for pseudo-steady irrotational flow.
 
From \eqref{eq: 2.7} it is easy to see that
 \begin{align}\label{eq: 2.8}
     \begin{cases}
     \bar{\partial}_\pm u=\cos(\sigma\pm\delta)+\dfrac{\cos \sigma}{\sin \delta}\bar{\partial}_\pm c+\dfrac{c \cos \alpha \bar{\partial}_\pm \beta-c \cos \beta \bar{\partial}_\pm \alpha}{2 \sin^2 \delta},\\
     \bar{\partial}_\pm v=\sin(\sigma\pm\delta)+\dfrac{\sin \sigma}{\sin \delta}\bar{\partial}_\pm c+\dfrac{c \sin \alpha \bar{\partial}_\pm \beta-c \sin \beta \bar{\partial}_\pm \alpha}{2 \sin^2 \delta}.
     \end{cases}
 \end{align}
Exploiting \eqref{eq: 2.8} in \eqref{eq: 2.6} yields
\begin{align}\label{eq: 2.9}
\begin{cases}
    \bar{\partial}_+ c=-\dfrac{\cos 2\delta}{\cot \delta}+\dfrac{c}{\sin 2 \delta}(\bar{\partial}_+\alpha-\cos 2\delta \bar{\partial}_+\beta),\vspace{0.2 cm}\\
    \bar{\partial}_- c=-\dfrac{\cos 2\delta}{\cot \delta}+\dfrac{c}{\sin 2 \delta}(\cos 2 \delta\bar{\partial}_-\alpha- \bar{\partial}_-\beta).
    \end{cases}
\end{align}
Now differentiating the pseudo-Bernoulli law \eqref{eq: 2.3} and using \eqref{eq: 2.8}, we have
\begin{align}\label{eq: 2.10}
    \left(\dfrac{1}{\sin^2\delta}+\kappa(\tau)\right)\bar{\partial}_\pm c=\cot \delta\left(\dfrac{c\bar{\partial}_\pm \delta}{\sin^2 \delta}-1\right),
\end{align}
where $\kappa(\tau)=\dfrac{-2p'(\tau)}{2p'(\tau)+\tau p''(\tau)}$.

Using \eqref{eq: 2.10} in \eqref{eq: 2.9} leads to the decompositions of the form
\begin{align}\label{eq: 2.11}
    \begin{cases}
    c\bar{\partial}_-\beta=\Upsilon(\tau, \delta)\cos^2 \delta( c\bar{\partial}_- \alpha- 2\sin^2 \delta),\\
    c\bar{\partial}_+\alpha=\Upsilon(\tau, \delta)\cos^2 \delta( c\bar{\partial}_+ \beta+ 2\sin^2 \delta),
    \end{cases}
\end{align}
where $\Upsilon(\tau, \delta)= m(\tau)-\tan^2 \delta= \dfrac{\kappa(\tau)-1}{\kappa(\tau)+1}-\tan^2 \delta$.
By exploiting \eqref{eq: 2.11} in \eqref{eq: 2.9} we can obtain
\begin{align}\label{eq: 2.12}
\begin{cases}
c \bar{\partial}_+\alpha&=-\dfrac{\Upsilon(\tau, \delta)}{2 \mu^2(\tau)}\sin{2 \delta}\bar{\partial}_+c= \dfrac{\tau^2p''(\tau)}{4c}\Upsilon \sin 2 \delta \bar{\partial}_+\tau,\\
c\bar{\partial}_+\beta&=-\dfrac{\tan \delta}{\mu^2(\tau)}\bar{\partial}_+c-2 \sin^2{\delta}=\dfrac{\tau^2p''(\tau)}{2c }\tan \delta \bar{\partial}_+\tau-2 \sin^2{\delta} ,\\
c\bar{\partial}_-\alpha&= \dfrac{\tan \delta}{\mu^2(\tau)}\bar{\partial}_-c+2 \sin^2{\delta}=-\dfrac{\tau^2p''(\tau)}{2c }\tan \delta \bar{\partial}_-\tau+2 \sin^2{\delta},\\
c \bar{\partial}_-\beta&=\dfrac{\Upsilon(\tau, \delta)}{2 \mu^2(\tau)}\sin{2 \delta}\bar{\partial}_-c= -\dfrac{\tau^2p''(\tau)}{4c}\Upsilon \sin 2 \delta \bar{\partial}_-\tau,
\end{cases}
\end{align} 
where $\mu^2(\tau)=\dfrac{1}{1+\kappa(\tau)}$.

Using \eqref{eq: 2.12} in \eqref{eq: 2.8}, we obtain the decompositions of velocity
\begin{align}\label{eq: 2.13}
    \begin{cases}
    \bar{\partial}_\pm u=\mp \dfrac{c}{\tau}\sin (\sigma\mp \delta)\bar{\partial}_\pm \tau, \vspace{0.1 cm}\\
    \bar{\partial}_\pm v=\pm \dfrac{c}{\tau}\cos (\sigma\mp \delta)\bar{\partial}_\pm \tau.
    \end{cases}
\end{align}
\subsection{Second order characteristic decompositions}
In this subsection, we derive characteristic decomposition forms for the variable $\rho$ which are very important to develop \textit{a priori} gradient estimates of solution in the interaction domain. We first use the following second-order normalized commutator relation from \cite{li2011characteristic}.
 \begin{p1}\label{p-2.1}
(Normalized commutator relation) \textit{ For any $I(\xi, \eta)$, we have}
\begin{align}\label{eq: 2.14}
    \bar{\partial}_-\bar{\partial}_+I-\bar{\partial}_+\bar{\partial}_-I=\dfrac{1}{\sin 2 \delta}\big[(\cos 2 \delta \bar{\partial}_+\beta-\bar{\partial}_-\alpha)\bar{\partial}_-I-(\bar{\partial}_+\beta-\cos 2 \delta\bar{\partial}_-\alpha)\bar{\partial}_+I\big].
\end{align}
 \end{p1}
\begin{p1}\label{p1}
For the variable $\rho$, we have the following second-order characteristic decompositions
\begin{align}\label{eq: 2.15}
    \begin{cases}
    c\bar{\partial}_+\bar{\partial}_-\rho=\bar{\partial}_-\rho\bigg[\sin 2 \delta +\dfrac{\tau^4 p''(\tau)}{4c \cos^2 \delta}\left(\bar{\partial}_-\rho+(f-1)\bar{\partial}_+\rho\right)\bigg], \vspace{0.2 cm}\\
    c\bar{\partial}_-\bar{\partial}_+\rho=\bar{\partial}_+\rho\bigg[\sin 2 \delta +\dfrac{\tau^4 p''(\tau)}{4c \cos^2 \delta}\left(\bar{\partial}_+\rho+(f-1)\bar{\partial}_-\rho\right)\bigg],
    \end{cases}
\end{align}
where 
\begin{align}
    f=2\sin^2 \delta-\dfrac{8p'(\tau)\cos^4 \delta}{\tau p''(\tau)}>0 ~as ~p''(\tau)>0~and~p'(\tau)<0.
\end{align}
\end{p1}
\begin{proof}
This decomposition was first derived in \cite{lai2019global} (see also \cite{lai2021interactions}). Here we only sketch the proof of this Proposition. We refer the reader to \cite{lai2019global} for more details.

We use the commutator relation \eqref{eq: 2.14} on $u$ and use \eqref{eq: 2.13} to obtain
\begin{align*}
    (\sin \alpha&+\sin \beta)\dfrac{1}{c\tau}\dfrac{d(c\tau)}{d\rho}\bar{\partial}_+\rho\bar{\partial}_-\rho+\sin\alpha\bar{\partial}_+\bar{\partial}_-\rho+\sin\beta\bar{\partial}_-\bar{\partial}_+\rho \\
    &=\dfrac{1}{\sin2\delta}\bigg[(\sin\alpha\bar{\partial}_-\alpha
    -\sin\alpha\cos2\delta\bar{\partial}_+\beta-\cos\alpha\sin2\delta\bar{\partial}_+\alpha)\bar{\partial}_-\rho \\
    &\hspace{1 cm}-(\sin\beta\bar{\partial}_+\beta-\sin\beta\cos2\delta\bar{\partial}_-\alpha+\cos\beta\sin2\delta\bar{\partial}_-\beta)\bar{\partial}_+\rho\bigg].
\end{align*}
Then one can apply the commutator relation on $\rho$ and use the relations \eqref{eq: 2.12} to prove this Proposition.
\end{proof}
\begin{p1}
If $\mathcal{F}(\rho)$ is any given smooth function then we have the following characteristic decompositions
\begin{align}
    \begin{cases}
    c\bar{\partial}_+\left(\dfrac{\mathcal{F}(\rho)\bar{\partial}_-\rho}{\sin^2 \delta}\right)=\dfrac{\tau p''(\tau)\bar{\partial}_-\rho}{4c \cos^2 \delta}\bigg[\dfrac{\mathcal{F}(\rho)\bar{\partial}_-\rho}{\sin^2 \delta}+\mathcal{G}\dfrac{\mathcal{F}(\rho)\bar{\partial}_+\rho}{\sin^2 \delta}\bigg],\\
    c\bar{\partial}_-\left(\dfrac{\mathcal{F}(\rho)\bar{\partial}_+\rho}{\sin^2 \delta}\right)=\dfrac{\tau p''(\tau)\bar{\partial}_+\rho}{4c \cos^2 \delta}\bigg[\dfrac{\mathcal{F}(\rho)\bar{\partial}_+\rho}{\sin^2 \delta}+\mathcal{G}\dfrac{\mathcal{F}(\rho)\bar{\partial}_-\rho}{\sin^2 \delta}\bigg],
    \end{cases}
\end{align}
where 
\begin{align}
    \mathcal{G}=2\sin^2 \delta-\dfrac{8p'(\tau)\cos^4 \delta}{\tau p''(\tau)}+\dfrac{c \mathcal{F}'(\rho)}{\mathcal{F}(\rho)}\dfrac{4c \cos^2 \delta}{\tau^4 p''(\tau)}-2 \cos^2 \delta+2\Upsilon \cos^4 \delta-1.
\end{align}
\end{p1}
\begin{proof}
The proof of this Proposition can be obtained by an easy manipulation on Proposition \ref{p1}, so we omit the details. 
\end{proof}
 \section{Interaction of planar and centered rarefaction waves}\label{3}
We solve the Riemann problem \eqref{eq: 1}-\eqref{eq: 2} by utilizing the characteristic decompositions \eqref{eq: 2.12} and \eqref{eq: 2.15} in this section. The expansion of gas through a bent corner can be considered as the interaction of a two-dimensional centered rarefaction wave with a planar rarefaction wave. For this reason, we first define these waves for 2-D Euler equations.
\subsection{Two-dimensional planar rarefaction wave}
Let us consider the system \eqref{eq: 1} with the initial data 
\begin{align}\label{eq: 3.1}
    (u, v, \tau)(x, y, 0)=\begin{cases}
    (u_0, 0, \tau_0),\hspace{1.5 cm}~ \mu x+\nu y<0,\\
    \mathrm{vacuum},\hspace{1.75 cm}~ \mu x+\nu y>0,
    \end{cases}
\end{align}
where $\mu^2+\nu^2=1$.

Let us use a transformation of coordinates of the form $\hat{x}=\mu x+\nu y,~\hat{y}=-\nu x+\mu y$. Accordingly, we denote $\hat{u}=\mu u+\nu v,~\hat{y}=-\nu u+\mu v$. Then one can solve a 1-D Riemann problem to obtain the planar rarefaction wave solution of the form
\begin{align}
    (\hat{u}, \hat{v}, \tau)(\hat{x}, \hat{y}, t)=\begin{cases}
    (u_0, 0, \tau_0),\hspace{1.75 cm} \hat{\xi}<\hat{\xi}_2,\\
    (u_r, 0, \tau_r)(\hat{\xi}),\hspace{1.25 cm}\hat{\xi}_2\leq \hat{\xi}\leq \hat{\xi}_1,\\
         \mathrm{vacuum},\hspace{2 cm}\hat{\xi}>\hat{\xi}_1,
    \end{cases}
\end{align}
where $\hat{\xi}=\hat{x}/t,~ \hat{\xi}_1=\underset{\tau_r\longrightarrow \infty}\lim u_r(\tau_r)$ and $\hat{\xi}_2=u_0-c(\tau_0)$. The functions $u_r(\hat{\xi})$ and $\tau_r(\hat{\xi})$ can be implicitly determined by
\begin{align}
    u_r=u_0-\displaystyle{\int_{\tau_0}^{\tau_r}{\dfrac{c}{\tau}}d\tau}~(\tau_0<\tau_r),~~\hat{\xi}=u_r-c(\tau_r).
\end{align}
Therefore, solution of the problem \eqref{eq: 1} and \eqref{eq: 3.1} is 
\begin{align}
    (u, v, \tau)(x, y, t)=\begin{cases}
    (u_0, 0, \tau_0),\hspace{1.75 cm}~\hat{\xi}<\hat{\xi}_2,\\
    (\mu u_r, \nu u_r, \tau_r)(\hat{\xi}),\hspace{0.55 cm}~~~\hat{\xi}_2\leq \hat{\xi}\leq \hat{\xi}_1,\\
    \mathrm{vacuum},\hspace{2 cm}~\hat{\xi}>\hat{\xi}_1. 
    \end{cases}
\end{align}
\subsection{Centered rarefaction wave}
Motivated by the steady flow in a sharp corner we are trying to construct the centered rarefaction wave for system \eqref{eq: 2.3} and \eqref{eq: 2.6} in this subsection. For simplicity, we assume that the wall angle $\theta<\alpha_v$ and $\alpha_v>2\bar{\delta}(\tau_0)+\dfrac{\pi}{2}$ such that $\bar{\delta}(\tau_0)=\arctan m(\tau_0)$, which means that cavitation will appear before the flow arrives the wall $OB$. Then, we define the $C_+$ type-centered rarefaction wave as in \cite{sheng2018interaction} as follows.

\begin{df}
Let $\Psi(t)$ be an angular domain of the form
\begin{align}
    \Psi(t)=\{(\xi, \eta)| \xi\in[0, t] , \xi \tan \alpha_v\leq \eta\leq \xi \tan \alpha_0, [\alpha_v, \alpha_0]\subset\left(-\dfrac{\pi}{2}, \dfrac{\pi}{2}\right)\}.
\end{align}
Then a function $(u, v, \tau)(\xi, \eta)$ is called a $C_+$ type centered rarefaction wave solution of the system \eqref{eq: 2.3} and \eqref{eq: 2.6} with $O(0, 0)$ as the center point if it satisfies the following properties:
\begin{enumerate}
    \item $(u, v, \tau)$ can be determined implicitly from the continuously differentiable functions $\eta=g(\xi, \alpha):=\xi \tan \alpha$ and $(u, v, \tau)(\xi, \eta)=(\hat{u}, \hat{v}, \hat{\tau})(\xi, \alpha)$ defined on a rectangular domain $\mathcal{D}(t)=\{(\xi, \alpha)|\xi\in [0, t], \alpha_v\leq \alpha \leq \alpha_0\}$. Further, for any $(\xi, \alpha)\in \Upsilon(t)/ \{\xi=0\}, g_\alpha(\xi, \alpha)>0$ holds.
    \item Function $(u, v, \tau)$ satisfies the system \eqref{eq: 2.3} and \eqref{eq: 2.6} on $\Psi(t)/ \{(0, 0)\}$.
    \item $\alpha= \alpha_0$ and $\alpha=\alpha_v$ are the $C_+$ characteristic angles when the flow arrives at the point $P$ and at the vacuum state, respectively and correspond to $\eta=\xi \tan \alpha_0$ and $\eta=\xi \tan \alpha_v$.
\end{enumerate}
\end{df}
\begin{figure}
    \centering
    \includegraphics[width=4 in]{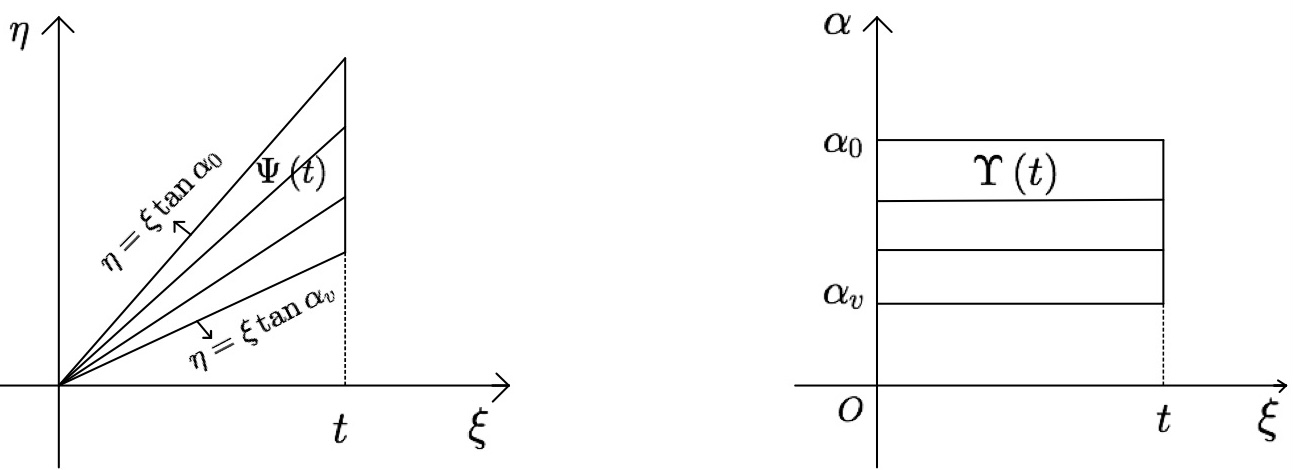}
    \caption{A $C_+$ type centered rarefaction wave in the $(\xi, \eta)$ and $(\xi, \alpha)$ plane}
    \label{fig: 2}
\end{figure}
The function $(\hat{u}, \hat{v}, \hat{\tau}, \hat{\phi})(0, \alpha):=(\hat{u}, \hat{v}, \hat{\tau}, \hat{\phi})( \alpha)$ is called the principal part of this $C_+$ type centered wave and $\alpha_0-\alpha_v$ is called amplitude of the centered wave. Then we have the following Lemma.
\begin{l1}
If $(u, v, \tau)(\xi, \eta)=(\hat{u}, \hat{v}, \hat{\tau})(\xi, \alpha),~\eta=\xi \tan \alpha, [\alpha_v, \alpha_0]\subset \left(-\dfrac{\pi}{2},\dfrac{\pi}{2}\right) $ is the $C_+$ type centered rarefaction wave solution of the system \eqref{eq: 2.3} and \eqref{eq: 2.6} in pseudo-supersonic domain, then the principal part $(\hat{u}, \hat{v}, \hat{\tau}, \hat{\phi})( \alpha)$ satisfy
\begin{align}
    \begin{cases}
    \dfrac{d\hat{u}}{d\alpha}+\tan \alpha \dfrac{d\hat{v}}{d\alpha}=0,~~ \hat{\phi}(\alpha)=const.,\\
    \dfrac{1}{2}(\hat{u}^2(\alpha)+\hat{v}^2(\alpha))+\displaystyle{\int_{\tau_0}^{\hat{\tau}}}\tau p'(\tau)d\tau=const.,\\
    \tan \alpha=\dfrac{\hat{u}(\alpha)\hat{v}(\alpha)+\hat{c}(\alpha)\sqrt{\hat{u}^2(\alpha)+\hat{v}^2(\alpha)-\hat{c}^2(\alpha)}}{\hat{u}^2(\alpha)-\hat{c}^2(\alpha)}.
    \end{cases}
\end{align}
\end{l1}
\begin{proof}
We substitute $(u, v, \tau)(\xi, \eta)=(\hat{u}, \hat{v}, \hat{\tau})(\xi, \alpha)$ and $\eta=\xi \tan \alpha$ into  pseudo-Bernoulli's law \eqref{eq: 2.3} and the system \eqref{eq: 2.6} to obtain 
\begin{align}\label{eq: 3.7}
    \begin{cases}
    \dfrac{\partial \hat{u}}{\partial \xi}+\tan \beta \dfrac{\partial \hat{v}}{\partial \xi}=0,\\
    \xi\sec^2 \alpha\dfrac{\partial \hat{v}}{\partial \xi}+\dfrac{\partial \hat{u}}{\partial \alpha}+\tan \alpha\dfrac{\partial \hat{v}}{\partial \alpha}=0,\\
    \dfrac{1}{2}((\hat{u}-\xi)^2+(\hat{v}-\eta)^2)+\displaystyle{\int_{\tau_0}^{\hat{\tau}}}\tau p'(\tau)d\tau +\hat{\phi}=const.
    \end{cases}
\end{align}
Now for the potential function $\phi$, we have $\phi_\xi=U$ and $\phi_\eta=V$. Then we have 
\begin{align}\label{eq: 3.8}
\xi\sec^2\alpha\dfrac{\partial \hat{\phi}}{\partial \xi}-\left(\tan \alpha\dfrac{\partial \hat{\phi}}{\partial \alpha}+\cot \sigma \dfrac{\partial \hat{\phi}}{\partial \alpha}\right)=0.
\end{align}
In view of \eqref{eq: 3.7}, \eqref{eq: 3.8} and considering $\xi\longrightarrow 0$, it is easy to obtain the required relations of the Lemma.
\end{proof}
In a similar manner, $C_-$ type centered wave solution can be also defined.
\subsection{Interaction of rarefaction waves and existence of local solution}
It is worth noting that when the gas expands into a vacuum from the sharp corner, the planar rarefaction wave $R_p$ and the centered rarefaction wave $R_c$ start interacting from the point $P$ in the self-similar plane where $P=\left(u_0-c_0, c_0\sqrt{\dfrac{u_0-c_0}{u_0+c_0}}\right)$. Therefore, we draw the $C_+$ wave characteristic curve $\overline{PQ}$ of $R_p$ and $C_-$ characteristic curve $\overline{PR}$ of $R_c$ to denote the interaction region $\Omega$ in $(\xi, \eta)$ plane which is bounded by curves $\overline{PQ}$, $\overline{PR}$ and vacuum boundary $\overline{QR}$.
\begin{figure}
    \centering
    \includegraphics[width= 4 in]{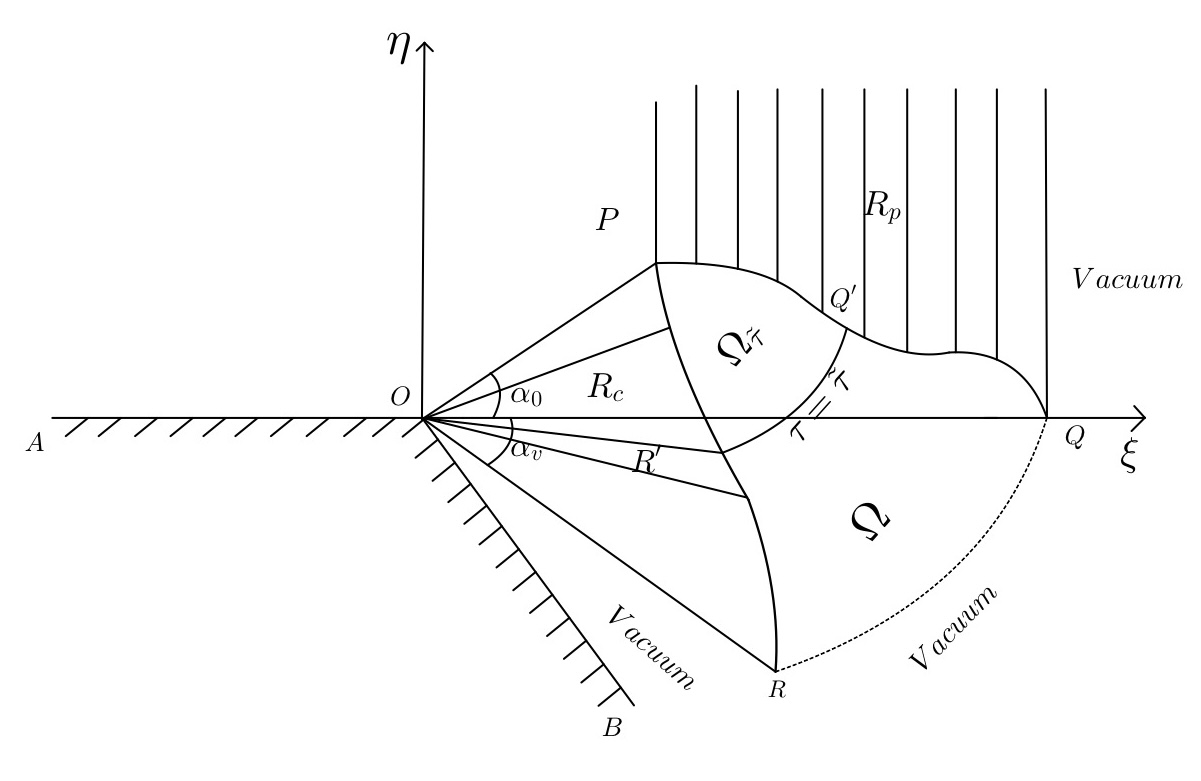}
    \caption{Interaction of planar rarefaction wave $R_p$ and centered rarefaction wave $R_c$}
    \label{fig: 3}
\end{figure}

The $C_+$ wave characteristic curve $\overline{PQ}$ can be expressed in the following form
\begin{align}
    \begin{cases}
    \xi=u_0-c-\displaystyle{\int_{\tau_0}^{\tau} \dfrac{c}{\tau}d\tau},\\
    \eta=\bigg\{\dfrac{\tau c}{\tau_0 c_0}\bigg[\left(\dfrac{2u_0c_0^2}{u_0+c_0}\right)-\tau_0c_0\left(\dfrac{c}{\tau}+\displaystyle{\int_{\tau_0}^{\tau}\dfrac{2c}{\tau^2}d\tau}\right)\bigg]\bigg\}^{1/2}.
    \end{cases}
\end{align}
From the expression of $\overline{PQ}$, it is easy to see that the concavity or convexity of the curve $\overline{PQ}$ depends on the choice of the sign of $m'(\tau)$ or $\kappa'(\tau)$ and can change accordingly for a general equation of state. Also, the point $Q$ is located at the $\xi$ axis when the gas enters into the vacuum and can be given by $Q(\xi, \eta)=\left(u_0+\displaystyle{\int_{0}^{\tau_0}\dfrac{c}{\tau}d\tau}, 0\right)$.

The values of $(u, v, \tau)$ on the $C_-$ characetristic curve $\overline{PR}$ and the $C_+$ characteristic curve $\overline{PQ}$ can be obtained by $R_c$ and $R_p$, respectively. Then the boundary data on $\overline{PQ}$ and $\overline{PR}$ is of the form
\begin{align}\label{eq: 3.10}
    (u, v, \tau)(\xi, \eta)=\begin{cases}
    (u_+, v_+, \tau_+)(\xi, \eta), ~~~\mathrm{on}~ \overline{PQ},\\
     (u_-, v_-, \tau_-)(\xi, \eta),~~~~\mathrm{on}~ \overline{PR}.
    \end{cases}
\end{align}
Noting that $\lambda_-= \dfrac{V^2-c^2}{UV+\sqrt{U^2+V^2-c^2}}$, we have $\lambda_-=-\infty$ or $\beta=-\dfrac{\pi}{2}$ on the curve $\overline{PQ}$. Then from \eqref{eq: 2.12} it is easy to see that  
\begin{align}\label{eq: 3.11}
   \bar{\partial}_+\beta=0, \bar{\partial}_+\rho=-\dfrac{2\sin 2 \omega c\rho^4}{p''(\tau)}<0
\end{align}
whenever $p''(\tau)>0$, which implies $\bar{\partial}_+\alpha>0$ along $\overline{PQ}$.

Similarly, across negative characteristic curve $\overline{PR}$, we must have
\begin{align}\label{eq: 3.12}
     \bar{\partial}_-\rho<0, \bar{\partial}_-\alpha<0.
\end{align}
Then we have the following Lemma:
\begin{l1}
(Local solution) \textit{For any $\tau\in  (\tau_0, \infty)$, the Goursat problem \eqref{eq: 2.1} and \eqref{eq: 3.10} admits a unique $C^1$ solution in the triangular domain $\Omega_{\tilde{\tau}}$ bounded by the $C_-$ characteristic curve $\overline{PR'}$, $C_+$ characteristic curve $\overline{PQ'}$ and the level curve $\tau=\tilde{\tau}$ connecting $Q'$ and $R'$ provided $\tilde{\tau}$ is sufficiently small. Moreover, this solution satisfies
}
\begin{align}
    \bar{\partial}_\pm\tau>0.
\end{align}
\end{l1}
\begin{proof}
From \cite{li1985boundary}, we know that for sufficiently small $\tau=\tilde{\tau}$, the Goursat problem \eqref{eq: 2.1} and \eqref{eq: 3.10} admits a unique $C^1$ solution in the domain closed by $\overline{PQ'}$, $\overline{PR'}$ and the level curve $\tau=\tilde{\tau}$. Further, using the boundary data \eqref{eq: 3.11} and \eqref{eq: 3.12} and the characteristic decomposition \eqref{eq: 2.15} we have $\bar{\partial}_\pm\rho<0$ or equivalently $\bar{\partial}_\pm\tau>0$.
\end{proof}
\section{Hyperbolicity and a priori $C^0$ and $C^1$ norm estimates}\label{4}
One of the main difficulties in the present problem is to control the hyperbolicity in the domain of determinacy while extending the local solution to the whole interaction region. In order to construct invariant regions of characteristic angles we use the ideas proposed in the work of Lai \cite{lai2021interactions} in this section.

Let us denote $\bar{\delta}(\tau)=\tan^{-1}\sqrt{m(\tau)}$, $\underset{\tau\rightarrow\infty}{\lim}\bar{\delta}(\tau)=\bar{\delta}^*$,  $\psi(\tau)=\bar{\delta}(\tau)-\bar{\delta}(\tau_0)+\displaystyle{\int_{\tau_0}^{\tau}}|\bar{\delta}'(\tau)|d\tau$,\\ $\chi(\tau)=\bar{\delta}(\tau)-\bar{\delta}(\tau_0)-\displaystyle{\int_{\tau_0}^{\tau}}|\bar{\delta}'(\tau)|d\tau$ and assume that 
\begin{align}\label{eq: 4.1}
    \begin{cases}
        m(\tau)>0, p'(\tau)<0, p''(\tau)>0,~~\mathrm{as}~\tau>\tau_0,\\
        2\bar{\delta}(\tau_0)+\chi(\tau)<\alpha_0+\dfrac{\pi}{2}<4\bar{\delta}(\tau_0).
    \end{cases}
\end{align}

Then we have the following Lemma:
\begin{l1}\label{l-4.1}
\textit{If the Goursat problem $\eqref{eq: 2.1}$, $\eqref{eq: 3.10}$ admits a $C^1$ solution in the domain $\Omega_{\tilde{\tau}}$ under the hypothesis \eqref{eq: 4.1} and $\tau=\tau_0, ~\tau=\tau_1$ and $\tau=\tau_2$ are the only points of local extrema of the curve $\bar{\delta}(\tau)$ in $\tau\in [\tau_0, \tau_2]$ such that $\tau_2<\tilde{\tau}\in(\tau_0, \infty)$ with $\bar{\delta}(\tau_0)<\bar{\delta}(\tau_1)$. Then there exists a positive constant $\epsilon_1$ such that \\ 
1). For any $\tau\in (\tau_0, \tau_1)$ if $m'(\tau)>0$ then we have
\begin{align}\label{eq: 4.2}
\begin{cases}
 \alpha\in \left(-\dfrac{\pi}{2}-\epsilon_1+2\bar{\delta}(\tau_0), \alpha_0+\epsilon_1+2(\bar{\delta}(\tau)-\bar{\delta}(\tau_0))\right)\vspace{0.2 cm}\\
 \beta\in \left(-\dfrac{\pi}{2}-\epsilon_1-2(\bar{\delta}(\tau)-\bar{\delta}(\tau_0)), \alpha_0+\epsilon_1-2\bar{\delta}(\tau_0)\right)
 \end{cases}
    \end{align}
for all $\tau\in[\tau_0, \tau_1]$ on $C_+$ characteristic $\overline{PQ_1}$ and $C_-$ characteristic $\overline{PR_1}$ where $Q_1$ and $R_1$ are the points corresponding to $\tau=\tau_1$ on the $C_+$ and $C_-$ characteristics $\overline{PQ}$ and $\overline{PR}$, respectively.\\
2). For any $\tau\in (\tau_1, \tau_2)$ if $m'(\tau)<0$ then we have
\begin{align}\label{eq: 4.3}
\begin{cases}
  \alpha\in \left(-\dfrac{\pi}{2}-\epsilon_1+2\bar{\delta}(\tau_0)+2(\bar{\delta
}(\tau)-\bar{\delta}(\tau_1)), \alpha_0+\epsilon_1\right)\vspace{0.2 cm}\\
\beta\in \left(-\dfrac{\pi}{2}-\epsilon_1, \alpha_0+\epsilon_1-2\bar{\delta}(\tau_0)-2(\bar{\delta}(\tau)-\bar{\delta}(\tau_1))\right)
\end{cases}
    \end{align}
for all $\tau\in(\tau_1, \tau_2]$ on $C_+$ characteristic $\overline{Q_1Q_2}$ and $C_-$ characteristic $\overline{R_1R_2}$ where $Q_2$ is a point on $\overline{PQ}$ while $R_2$ is a point on $\overline{PR}$ corresponding to $\tau=\tau_2$, respectively.}
\end{l1}
\begin{proof}
Case 1. Let us assume that for $\tau\in (\tau_0, \tau_1), m'(\tau)>0$ or equivalently $\bar{\delta}'(\tau)>0$. Then using the fact that $\bar{\partial}_\pm \tau>0$ in $\Omega_{\tilde{\tau}}$ we have $\bar{\partial}_\pm \bar{\delta}>0$ for all $\tau\in (\tau_0, \tau_1)$. Also for $\bar{\delta}'(\tau)>0$, we have $\chi(\tau)=0$. Then on the $C_+$ characteristic $\overline{PQ_1}$, we must have  
\begin{align*}
(\alpha, \beta)(P)=\left(\alpha_0, -\dfrac{\pi}{2}\right)\subset &\left(-\dfrac{\pi}{2}-\epsilon_1+2\bar{\delta}(\tau_0), \alpha_0+\epsilon_1+2(\bar{\delta}(\tau)-\bar{\delta}(\tau_0))\right)\\
&\times \left(-\dfrac{\pi}{2}-\epsilon_1-2(\bar{\delta}(\tau)-\bar{\delta}(\tau_0)), \alpha_0+\epsilon_1-2\bar{\delta}(\tau_0)\right)
\end{align*}
for a sufficiently small positive constant $\epsilon_1$. 
\begin{figure}
    \centering
    \includegraphics[width= 2.3 in]{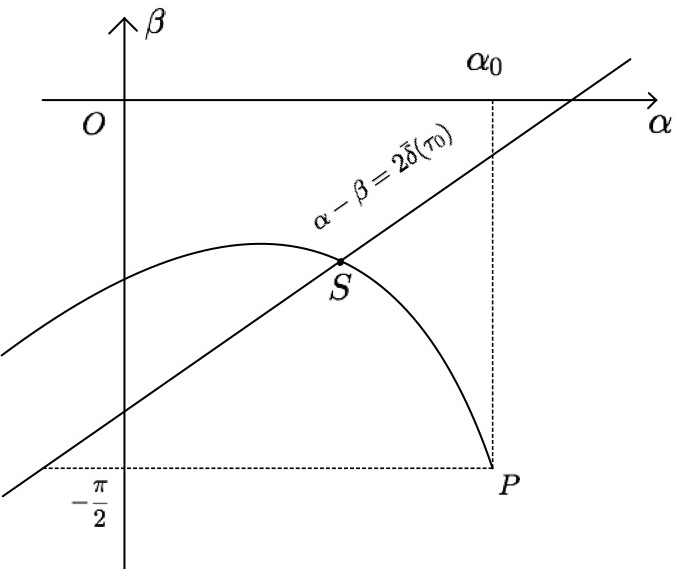}
    \caption{Beta curve $\beta=\beta(\alpha)$}
    \label{fig: 5}
\end{figure}

If there exists a point $I$ on $\overline{PQ_1}$ such that $\alpha(I)=-\dfrac{\pi}{2}-\epsilon_1+2\bar{\delta}(\tau_0)$. Then we have 
$$\delta(I)=\dfrac{\alpha-\beta}{2}=\dfrac{2\bar{\delta}(\tau_0)-\epsilon_1}{2}<\bar{\delta}(\tau_0)\leq \bar{\delta}(\tau).$$
Hence, we use \eqref{eq: 2.12} to obtain
$$c\bar{\partial}_+\alpha(I)=\dfrac{\tau^2p''(\tau)}{4c}(\tan^2 \bar{\delta}(\tau)-\tan^2 \delta)\sin 2 \delta \bar{\partial}_+\tau>0.$$
Again, if there exists a point $J$ on $\overline{PQ_1}$ such that $\alpha(J)=\alpha_0+\epsilon_1+2(\bar{\delta}(\tau)-\bar{\delta}(\tau_0))$ then using the hypothesis \eqref{eq: 4.1} we have
$$\delta(J)=\dfrac{\alpha_0+\epsilon_1+2(\bar{\delta}(\tau)-\bar{\delta}(\tau_0))+\pi/2}{2}>\bar{\delta}(\tau),$$
which means that 
$$\bar{\partial}_+\alpha(J)<0.$$ 
Combining the above results it is easy to see that $ -\dfrac{\pi}{2}-\epsilon_1+2\bar{\delta}(\tau_0)<\alpha< \alpha_0+\epsilon_1+2(\bar{\delta}(\tau)-\bar{\delta}(\tau_0))$ for all $\tau\in [\tau_0, \tau_1]$. Hence we have proved the first part of Lemma for positive characteristic $\overline{PQ_1}$.

Now for $C_-$ characteristic curve $\overline{PR_1}$, we use the first order decomposition of $\beta$ 
\begin{align}
    c\bar{\partial}_-\beta= -\dfrac{\tau^2p''(\tau)}{4c}(\tan^2 \bar{\delta}(\tau)-\tan^2 \delta)\sin 2 \delta \bar{\partial}_-\tau.
\end{align}
Then in view of $\delta(\tau_0)>\bar{\delta}(\tau_0)$ and $\bar{\partial}_-\tau>0$, we have $\bar{\partial}_-\beta>0$ at the point $P$ which together with the fact that $\bar{\partial}_-\alpha<0$ implies that $\bar{\partial}_-\delta<0$ at the point $P$ or in other words $\beta'(\alpha)|_P=\dfrac{\bar{\partial}_-\beta}{\bar{\partial}_-\alpha}\bigg|_P<0$. 

Now we consider the properties of the curve $\beta=\beta(\alpha)$ along $\overline{PR_1}$. Let $l$ be the line $\alpha-\beta=2\bar{\delta}(\tau_0)$ in $(\alpha, \beta)$ plane. It is easy to see that below the line $l$ we have $\delta(\tau)>\bar{\delta}(\tau)$ therefore $\bar{\partial}_-\beta>0$ or $\beta'(\alpha)<0$; see Figure \ref{fig: 5}. Further, if $S$ is the point of intersection of the two curves $\beta(\alpha)$ and $l$, then $\beta'(\alpha)|_S=0$. It is worth noting from the fact $\bar{\partial}_-\alpha<0$ that the two curves have only one intersection point. Moreover, $\beta=\beta(\alpha)$ must pass through $l$ if $(\alpha, \beta)$ is located above $l$. In the region above $l$, we must have $\delta(\tau)<\bar{\delta}(\tau)$ or $\beta'(\alpha)>0$. Now on the curve $\overline{PR_1}$, we have $\bar{\partial}_-\bar{\delta}(\tau)>0$ or $\bar{\delta}(\tau_0)<\bar{\delta}(\tau)$ for any $\tau\in [\tau_0, \tau_1]$, which together with the hypothesis \eqref{eq: 4.1} implies that $-\dfrac{\pi}{2}-\epsilon_1+2\bar{\delta}(\tau_0)<\alpha< \alpha_0+\epsilon_1+2(\bar{\delta}(\tau)-\bar{\delta}(\tau_0))$ and $-\dfrac{\pi}{2}-\epsilon_1-2(\bar{\delta}(\tau)-\bar{\delta}(\tau_0))<\beta< \alpha_0+\epsilon_1-2\bar{\delta}(\tau_0)$ on $\overline{PR_1}$.\\\\
Case 2. Now let us consider the case $m'(\tau)<0$ or $\bar{\delta}'(\tau)<0$ when $\tau\in (\tau_1, \tau_2)$. Then we are going to prove that 
\begin{align*}
\left(\alpha, \beta\right)\in \left(-\dfrac{\pi}{2}-\epsilon_1+2\bar{\delta}(\tau_0)+2(\bar{\delta
}(\tau)-\bar{\delta}(\tau_1)), \alpha_0+\epsilon_1\right)
\times \left(-\dfrac{\pi}{2}-\epsilon_1, \alpha_0+\epsilon_1-2\bar{\delta}(\tau_0)-2(\bar{\delta}(\tau)-\bar{\delta}(\tau_1))\right).
\end{align*}

Again on $C_+$ characteristic $\overline{Q_1Q_2}$ starting from the point $Q_1$, we have $\bar{\partial}_+\beta=0$ so that if there exists a point on $\overline{Q_1Q_2}$ such that $\alpha=-\dfrac{\pi}{2}-\epsilon_1+2\bar{\delta}(\tau_0)+2(\bar{\delta}(\tau)-\bar{\delta}(\tau_1))$, then  $$\delta=\dfrac{2\bar{\delta}(\tau_0)+2(\bar{\delta}(\tau)-\bar{\delta}(\tau_1))-\epsilon_1}{2}<\bar{\delta}(\tau)$$ as $\bar{\delta}(\tau_0)<\bar{\delta}(\tau_1)$
which implies that $\bar{\partial}_+\alpha>0.$ 

On a similar ground we can prove that $\bar{\partial}_+\alpha<0$ under the hypothesis \eqref{eq: 4.1} if there exists a point on $\overline{Q_1Q_2}$ such that $\alpha= \alpha_0+\epsilon_1$.

For the $C_-$ characteristics $\overline{R_1R_2}$ starting from the point $R_1$, we have $\beta'(\alpha)|_{R_1}>0$. Therefore, a similar argument as in case 1 of this Lemma is enough to prove this part as well so we don't repeat the arguments again for the sake of brevity.
\end{proof}
We can use the Lemma \ref{l-4.1} to obtain an invariant square for the characteristic angles in the following Proposition. 
\begin{figure}
    \centering
    \includegraphics[width=4.8 in]{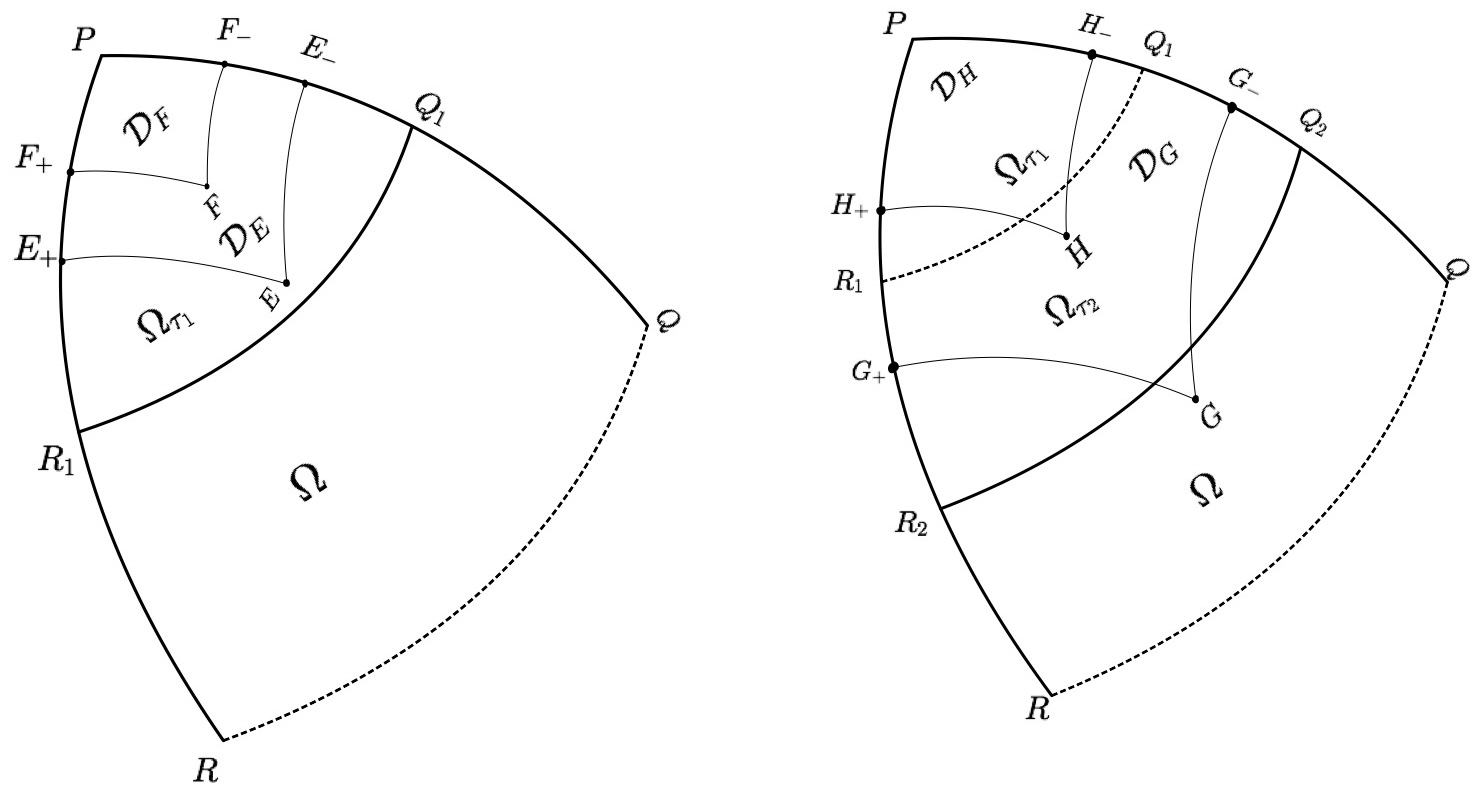}
    \caption{The regions $\Omega_{\tau_1}$ and $\Omega_{\tau_2}$.}
    \label{fig: 6}
\end{figure}
\begin{p1}\label{p-1}
\textit{Let $\tau=\tau_0, \tau=\tau_1$ and $\tau=\tau_2$ are the only points of local extrema of the curve $\bar{\delta}(\tau)$ in $\tau\in [\tau_0,\tau_2]$  such that $\tau_2<\tilde{\tau}\in(\tau_0, \infty)$ with $\bar{\delta}(\tau_0)<\bar{\delta}(\tau_1)$ then
\begin{enumerate}
    \item For any $\tau\in (\tau_0, \tau_1)$ if $m'(\tau)>0$ then 
    \begin{align}\label{eq: 4.5}
\begin{cases}
 \alpha\in \left(-\dfrac{\pi}{2}-\epsilon_1+2\bar{\delta}(\tau_0), \alpha_0+\epsilon_1+2(\bar{\delta}(\tau)-\bar{\delta}(\tau_0))\right)\vspace{0.2 cm}\\
 \beta\in \left(-\dfrac{\pi}{2}-\epsilon_1-2(\bar{\delta}(\tau)-\bar{\delta}(\tau_0)), \alpha_0+\epsilon_1-2\bar{\delta}(\tau_0)\right)
 \end{cases}
    \end{align}
for all $\tau\in[\tau_0, \tau_1]$ in $\Omega_{\tau_1}$ where $\Omega_{\tau_1}$ is a closed region bounded by the curves $\overline{PQ_1}$, $\overline{PR_1}$ and $\overline{Q_1R_1}$; see Figure \ref{6}.
    \item For any $\tau\in (\tau_1, \tau_2)$ if $m'(\tau)<0$ then 
\begin{align}\label{eq: 4.6}
\begin{cases}
  \alpha\in \left(-\dfrac{\pi}{2}-\epsilon_1+2\bar{\delta}(\tau_0)+2(\bar{\delta
}(\tau)-\bar{\delta}(\tau_1)), \alpha_0+\epsilon_1\right)\vspace{0.2 cm}\\
\beta\in \left(-\dfrac{\pi}{2}-\epsilon_1, \alpha_0+\epsilon_1-2\bar{\delta}(\tau_0)-2(\bar{\delta}(\tau)-\bar{\delta}(\tau_1))\right)
\end{cases}
    \end{align}
for all $\tau\in(\tau_1, \tau_2]$ in $\Omega_{\tau_2}/\Omega_{\tau_1}$ where $\Omega_{\tau_2}$ is a closed region bounded by the curves $\overline{PQ_2}$, $\overline{PR_2}$ and $\overline{Q_2R_2}$; see Figure \ref{6}.
\end{enumerate}
}
%1). Let E be an arbitrary point in $D_{\tau_1}$ and \\
%2). Let $G$ be an arbitrary point in $D_{\tau_2}/D_{\tau_1}$. The $C_+$ characteristic curve passing through $G$ intersects with $IH$ at the point $G_+$, and the $C_-$ characteristic curve passing through $G$ intersects with $IK$ at the point $G_-$. Let $\Delta_G$ be the closed domain bounded by the characteristic curves $IG_+$, $IG_-$, $G_+G$ and $G_-G$. Then if $(\alpha, \beta)(\xi, \eta)\in \Gamma_2(\tau)$ as $(\xi, \eta)\in \Delta_G/{G}$ then $(\alpha, \beta)\in \Gamma_2(\tau_G)$ for all $(\xi, \eta)\in \Delta_G$, where $\tau_1<\tau_G\leq \tau$.}
\end{p1}
\begin{proof}
Case 1. To prove the first part of Proposition we need to prove that for any arbitrary point $E$ in $\Omega_{\tau_1}, (\alpha, \beta)(E)\in \Gamma(\tau)$, where 
\begin{align*}
\Gamma: \left(-\dfrac{\pi}{2}-\epsilon_1+2\bar{\delta}(\tau_0), \alpha_0+\epsilon_1+2(\bar{\delta}(\tau)-\bar{\delta}(\tau_0))\right)
\times \left(-\dfrac{\pi}{2}-\epsilon_1-2(\bar{\delta}(\tau)-\bar{\delta}(\tau_0)), \alpha_0+\epsilon_1-2\bar{\delta}(\tau_0)\right).
\end{align*}

Let us assume that $C_+$ characteristic curve passing through $E$ intersects $\overline{PR_1}$ at $E_+$ while the $C_-$ characteristic curve passing through $E$ intersects $\overline{PQ_1}$ at $E_-$ and $\mathcal{D}_E$ is a closed domain bounded by characteristic curves $\overline{PE_+}$, $\overline{PE_-}$, $\overline{E_+E}$ and $\overline{E_-E}$; see Figure \ref{6}. Then we claim that if $(\alpha, \beta)\in \Gamma(\tau)$ as $(\xi, \eta)\in \mathcal{D}_E/\{E\}$ then $(\alpha, \beta)\in \Gamma(\tau_E)$ for all $(\xi, \eta)\in \mathcal{D}_E$, where $\tau_E\leq \tau$. 

\begin{figure}
    \centering
    \includegraphics[width= 3.8 in]{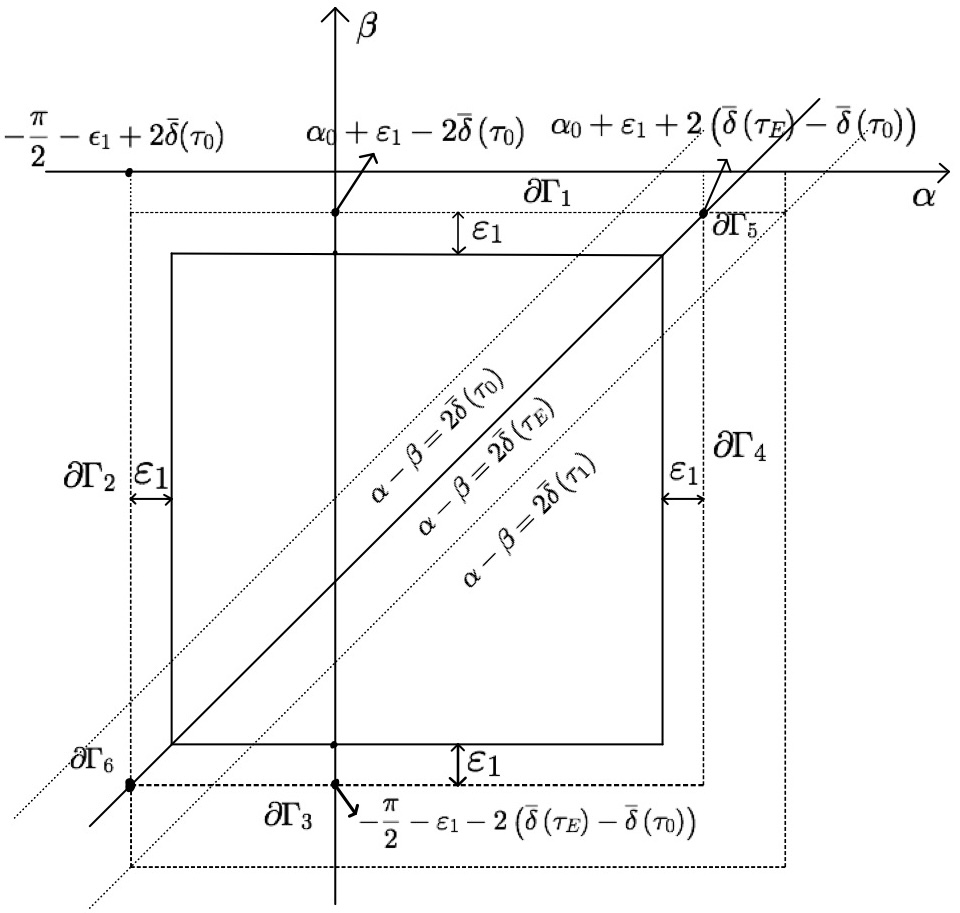}
    \caption{Invariant square for the case $m'(\tau)>0$}
    \label{fig: 6}
\end{figure}
We now divide the boundary of $\Gamma(\tau_E)$ into the following six parts (see Figure \ref{fig: 6}):
\begin{itemize}
    \item $\partial\Gamma_1=\{-\dfrac{\pi}{2}-\epsilon_1+2\bar{\delta}(\tau_0)\leq \alpha \leq \alpha_0+\epsilon_1+2(\bar{\delta}(\tau_E)-\bar{\delta}(\tau_0)), \beta= \alpha_0+\epsilon_1-2\bar{\delta}(\tau_0)\}$
    \item $\partial\Gamma_2=\{\alpha=-\dfrac{\pi}{2}-\epsilon_1+2\bar{\delta}(\tau_0), -\dfrac{\pi}{2}-\epsilon_1-2(\bar{\delta}(\tau_E)-\bar{\delta}(\tau_0)) \leq\beta\leq \alpha_0+\epsilon_1-2\bar{\delta}(\tau_0) \}$
    \item $\partial\Gamma_3=\{-\dfrac{\pi}{2}-\epsilon_1+2\bar{\delta}(\tau_0)\leq \alpha \leq \alpha_0+\epsilon_1+2(\bar{\delta}(\tau_E)-\bar{\delta}(\tau_0)), \beta=  -\dfrac{\pi}{2}-\epsilon_1-2(\bar{\delta}(\tau_E)-\bar{\delta}(\tau_0)) \}$
    \item $\partial\Gamma_4=\{\alpha=\alpha_0+\epsilon_1+2(\bar{\delta}(\tau_E)-\bar{\delta}(\tau_0)), -\dfrac{\pi}{2}-\epsilon_1-2(\bar{\delta}(\tau_E)-\bar{\delta}(\tau_0)) \leq\beta\leq \alpha_0+\epsilon_1-2\bar{\delta}(\tau_0) \}$
    \item $\partial\Gamma_5= \{\alpha=\alpha_0+\epsilon_1+2(\bar{\delta}(\tau_E)-\bar{\delta}(\tau_0)),\beta= \alpha_0+\epsilon_1-2\bar{\delta}(\tau_0) \}$
    \item $\partial\Gamma_6=\{\alpha=-\dfrac{\pi}{2}-\epsilon_1+2\bar{\delta}(\tau_0), \beta= -\dfrac{\pi}{2}-\epsilon_1-2(\bar{\delta}(\tau_E)-\bar{\delta}(\tau_0))\}$
\end{itemize}

From Lemma \ref{l-4.1} we have $(\alpha, \beta)\in \Gamma(\tau_E)$ on $\overline{PE_+}$ and $\overline{PE_-}$. Then if we assume that the Proposition is not true then there must exist a point $F\in \mathcal{D}_E$ such that $(\alpha, \beta)\in \Gamma(\tau_E)$ for all $(\xi, \eta)\in \mathcal{D}_F/\{F\}$ and $(\alpha, \beta)(F)\in \bigcup_{i=1}^6 \partial\Gamma_i$ where $\mathcal{D}_{F}$ can be defined in a same manner as $\mathcal{D}_{E}$.

Let us assume that $(\alpha, \beta)(F)\in \partial\Gamma_1$. Then using the hypothesis $(\alpha, \beta)\in\Gamma(\tau)$ as $(\xi, \eta)\in \mathcal{D}_{E}/\{E\}$ and noting that $\partial\Gamma_1\subset \partial \Gamma(\tau)$, we obtain $E=F$. Now since $\delta< \bar{\delta}(\tau_E)$ in $\partial\Gamma_1$ so we have
\begin{align*}
    c\bar{\partial}_-\beta(F)<0.
\end{align*}
However, from $(\alpha, \beta)\in \Gamma(\tau_E)$ as $(\xi, \eta)\in \mathcal{D}_F/\{F\}$ we must have $\bar{\partial}_-\beta(F)\geq 0$, which yields a contradiction. Hence $(\alpha, \beta)\not\in\partial\Gamma_1$. Similarly, one can prove that  $(\alpha, \beta)\not\in\partial\Gamma_2$.

Now let us assume that $(\alpha, \beta)(F)\in \partial\Gamma_3$. Then noting that $\delta>\bar{\delta}(E)\geq \bar{\delta}(F)$ on $\partial\Gamma_3$, we must have
\begin{align*}
    c\bar{\partial}_-\beta(F)>0.
\end{align*}
However, according to $(\alpha, \beta)\in \Gamma(\tau_E)$ as $(\xi, \eta)\in \mathcal{D}_F/\{F\}$ we must have $\bar{\partial}_-\beta(F)\leq 0$, which is a contradiction. Hence $(\alpha, \beta)\not\in\partial\Gamma_3$. Similarly, one can prove that  $(\alpha, \beta)\not\in\partial\Gamma_4$.

Now if $(\alpha, \beta)\in \partial\Gamma_5$, then we define a function $\hat{\alpha}$ on $\overline{F_+F}$ by the following relations
\begin{align}\label{eq: 4.9}
    \begin{cases}
    c\bar{\partial}_+\hat{\alpha}=\dfrac{\tau^2p''(\tau)}{4c}\bigg[\tan^2\bar{\delta}(\tau_E)-\tan^2\left(\dfrac{\hat{\alpha}-\alpha_0+2\bar{\delta}(\tau_0)-\epsilon_1}{2}\right)\bigg]\bar{\partial}_+\tau,\\
    \hat{\alpha}(F_+)=\alpha(F_+).
    \end{cases}
\end{align}
According to hypothesis that $(\alpha, \beta)\in \Gamma(\tau_E)$ as $(\xi, \eta)\in \mathcal{D}_F/\{F\}$, it is easy to see that $\alpha(F_+)=\hat{\alpha}(F_+)< \alpha_0+\epsilon+2(\bar{\delta}(\tau_E)-\bar{\delta}(\tau_0))$. Then
integrating \eqref{eq: 4.9} along the positive characteristic curve $\overline{F_+F}$ from $F_+$ to $F$, we obtain
$$\hat{\alpha}(F)<\alpha_0+\epsilon+2(\bar{\delta}(\tau_E)-\bar{\delta}(\tau_0)).$$
Combining \eqref{eq: 4.9} and \eqref{eq: 2.12}, we have
\begin{align}
    \begin{cases}
    c\bar{\partial}_+(\hat{\alpha}-\alpha)=\dfrac{\tau^2p''(\tau)}{4c}\big[\tan^2\bar{\delta}(\tau_E)-\tan^2\bar{\delta}(\tau)\big]\bar{\partial}_+\tau\\
    ~~~~~~~~~~~~~~-\dfrac{\tau^2p''(\tau)}{4c}\big[\tan^2\left(\dfrac{\hat{\alpha}-\alpha_0+2\bar{\delta}(\tau_0)-\epsilon_1}{2}\right)-\tan^2\left(\dfrac{\alpha-\beta}{2}\right)\big]\bar{\partial}_+\tau,\\
    (\hat{\alpha}-\alpha)(F_+)=0.
    \end{cases}
\end{align}
Now if $\alpha=\hat{\alpha}$ then from the assumption that $(\alpha, \beta)\in \Gamma(\tau_E)$ as $(\xi, \eta)\in \mathcal{D}_F/\{F\}$ we must have\\
$\bigg[\tan^2\left(\dfrac{\hat{\alpha}-\alpha_0+2\bar{\delta}(\tau_0)-\epsilon_1}{2}\right)-\tan^2\left(\dfrac{\alpha-\beta}{2}\right)\bigg]>0$ or 
$\bar{\partial}_+(\hat{\alpha}-\alpha)>0$.

Therefore, we have $\alpha(F)<\hat{\alpha}(F)<\alpha_0+\epsilon_1+2(\bar{\delta}(\tau_E)-\bar{\delta}(\tau_0))$, which leads to a contradiction. Hence, $(\alpha, \beta)\notin \partial\Gamma_5$. On a similar ground, one can prove that $(\alpha, \beta)\notin \partial\Gamma_6$.

Combining all these results it can be stated that there is no such point $F\in \mathcal{D}_E$, hence we have proved the first part of Proposition.\\\\
Case 2. Now we proceed to the proof of the second part of the Proposition. Let us assume that the curve $\bar{\delta}(\tau)$ be monotonically decreasing in $(\tau_1, \tau_2)$. Then we need to prove that $(\alpha, \beta)(G)\in \Gamma'(\tau)$ for any $G\in \Omega_{\tau_2}/\Omega_{\tau_1}$, where \begin{align*}
\Gamma'(\tau): \left(-\dfrac{\pi}{2}-\epsilon_1+2\bar{\delta}(\tau_0)+2(\bar{\delta
}(\tau)-\bar{\delta}(\tau_1)), \alpha_0+\epsilon_1\right)
\times \left(-\dfrac{\pi}{2}-\epsilon_1, \alpha_0+\epsilon_1-2\bar{\delta}(\tau_0)-2(\bar{\delta}(\tau)-\bar{\delta}(\tau_1))\right).
\end{align*}
Similar to first part we divide the boundary of $\Gamma'(\tau_G)$ into the following six parts
\begin{figure}
    \centering
    \includegraphics[width=3.8 in]{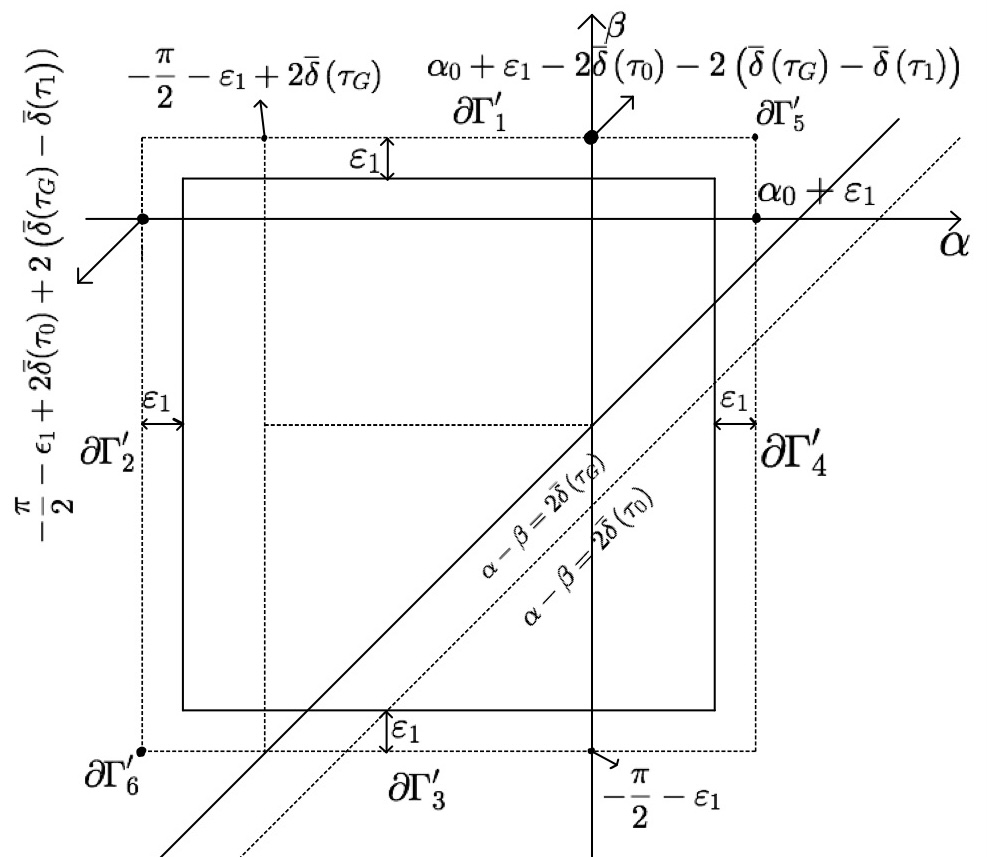}
    \caption{Invariant square for the case $m'(\tau)<0$}
    \label{fig:my_label}
\end{figure}
\begin{itemize}
    \item $\partial\Gamma_1'=\bigg\{-\dfrac{\pi}{2}-\epsilon_1+2\bar{\delta}(\tau_0)+2(\bar{\delta}(\tau_G)-\bar{\delta}(\tau_1))\leq \alpha \leq \alpha_0+\epsilon_1, \beta=\alpha_0+\epsilon_1-2\bar{\delta}(\tau_0)-2(\bar{\delta}(\tau_G)-\bar{\delta}(\tau_1))\bigg\}$;
    \item $\partial\Gamma_2'=\bigg\{\alpha=-\dfrac{\pi}{2}-\epsilon_1+2\bar{\delta}(\tau_0)+2(\bar{\delta}(\tau_G)-\bar{\delta}(\tau_1)),-\dfrac{\pi}{2}-\epsilon_1\leq \beta\leq \alpha_0+\epsilon_1-2\bar{\delta}(\tau_0)-2(\bar{\delta}(\tau_G)-\bar{\delta}(\tau_1))\bigg\}$;
     \item $\partial\Gamma_3'=\bigg\{-\dfrac{\pi}{2}-\epsilon_1+2\bar{\delta}(\tau_0)+2(\bar{\delta}(\tau_G)-\bar{\delta}(\tau_1))\leq \alpha \leq \alpha_0+\epsilon_1, \beta=-\dfrac{\pi}{2}-\epsilon_1\bigg\}$;
     \item $\partial\Gamma_4'=\bigg\{\alpha=\alpha_0+\epsilon_1,-\dfrac{\pi}{2}-\epsilon_1\leq \beta\leq \alpha_0+\epsilon_1-2\bar{\delta}(\tau_0)-2(\bar{\delta}(\tau_G)-\bar{\delta}(\tau_1))\bigg\}$;
     \item $\partial\Gamma_5'=\bigg\{\alpha=\alpha_0+\epsilon_1,\beta=\alpha_0+\epsilon_1-2\bar{\delta}(\tau_0)-2(\bar{\delta}(\tau_G)-\bar{\delta}(\tau_1))\bigg\}$;
     \item $\partial\Gamma_6'=\bigg\{\alpha=-\dfrac{\pi}{2}-\epsilon_1+2\bar{\delta}(\tau_0)+2(\bar{\delta}(\tau_G)-\bar{\delta}(\tau_1)),\beta=-\dfrac{\pi}{2}-\epsilon_1 \bigg\}$.
\end{itemize}
Similar to the first part of the proof let us assume that the Proposition is not true then there must exist a point $H\in \mathcal{D}_G$ such that $(\alpha, \beta)\in \Gamma'(\tau_G)$ for all $(\xi, \eta)\in \mathcal{D}_H/\{H\}$ and $(\alpha, \beta)(H)\in \bigcup_{i=1}^6\partial\Gamma_i'$ where $\mathcal{D}_G$ and $\mathcal{D}_H$ can be defined similarly to $\mathcal{D}_E$.

First let us suppose that $(\alpha, \beta)(H)\in \partial\Gamma_1'$. Then since $\delta(H)<\bar{\delta}(\tau_G)\leq\bar{\delta}(\tau_H)$, we have $\bar{\partial}_-\beta(H)<0$, which is a contradiction since $\bar{\partial}_-\beta(H)\geq0$ as $(\alpha, \beta)\in \Gamma'(\tau_G)$ for all $(\xi, \eta)\in \mathcal{D}_H/\{H\}$. Therefore, $(\alpha, \beta)(H)\notin \partial\Gamma_1'$. Similarly, one can prove that $(\alpha, \beta)(H)\notin \partial\Gamma_2'$.

Now let us suppose that $(\alpha, \beta)(H)\in \partial\Gamma_3'$. Then by the assumption of Proposition, we must have $G=H$, since $\partial\Gamma_3'\subset \partial \Gamma'(\tau)$. Hence we have $\bar{\delta}(\tau_H)=\bar{\delta}(\tau_G)<\delta(H)$, Consequently, we obtain that $\bar{\partial}_-\beta(H)>0$ which leads to a contradiction since $\bar{\partial}_-\beta(H)\leq0$ as $(\alpha, \beta)\in \Gamma'(\tau_G)$ for all $(\xi, \eta)\in \mathcal{D}_H/\{H\}$. Similarly, one can prove that $(\alpha, \beta)\notin \partial \Gamma_4'$.

Now let us assume that $(\alpha, \beta)\in \partial \Gamma_6'$, then we define the function $\hat{\alpha}$ which satisfies
\begin{align}\label{eq: 4.10}
    \begin{cases}
    c\bar{\partial}_+\hat{\alpha}=\dfrac{\tau^2p''(\tau)}{4c}\bigg[\tan^2\bar{\delta}(\tau_G)-\tan^2\left(\dfrac{\hat{\alpha}+\frac{\pi}{2}+\epsilon_1}{2}\right)\bigg]\bar{\partial}_+\tau,\\
    \hat{\alpha}(H_+)=\alpha(H_+),
    \end{cases}
\end{align}
on $\overline{H_+H}$. Then since $(\alpha,\beta)(H_+)\in\Gamma'(\tau_G)$, we have
\begin{align}
    -\dfrac{\pi}{2}-\epsilon_1+2\bar{\delta}(\tau_0)+2(\bar{\delta}(\tau)-\bar{\delta}(\tau_1))<\alpha(H_+)=\hat{\alpha}(H_+).
\end{align}
Then combining \eqref{eq: 2.12} and \eqref{eq: 4.10} we have
\begin{align}
    \begin{cases}
    c\bar{\partial}_+(\hat{\alpha}-\alpha)=\dfrac{\tau^2p''(\tau)}{4c}\bigg[\tan^2\bar{\delta}(\tau_G)-\tan^2\bar{\delta}(\tau)-\tan^2\left(\dfrac{\hat{\alpha}+\frac{\pi}{2}+\epsilon_1}{2}\right)+\tan^2\left(\dfrac{\alpha-\beta}{2}\right)\bigg]\bar{\partial}_+\tau,\\
    (\hat{\alpha}-\alpha)(H_+)=0.
    \end{cases}
\end{align}
Then noting that $\bar{\delta}(\tau)>\bar{\delta}(\tau_G)$ on $\overline{H_+H}$ and the hypothesis on $H$, we can integrate above to obtain
\begin{align}
    -\dfrac{\pi}{2}-\epsilon_1+2\bar{\delta}(\tau_0)+2(\bar{\delta}(\tau)-\bar{\delta}(\tau_1))<\alpha(H)<\hat{\alpha}(H),
\end{align}
which yields a contradiction, hence $(\alpha, \beta)\notin \partial \Gamma_6'$. In a similar manner, one can prove that $(\alpha, \beta)\notin \partial \Gamma_5'$. Hence proved.
\end{proof}
In view of the boundary data \eqref{eq: 3.11} and \eqref{eq: 3.12} we define 
\begin{align}
  M_1(\tau)=\min\{\bar{\partial}_+\rho|_{\overline{PQ}}, \bar{\partial}_-\rho|_{\overline{PR}}\} ~~\mathrm{and}~M_2(\tau)=\dfrac{\rho^nM_1(\tau)}{\sin^2\delta}.
\end{align}
Then we have the following Lemma.
\begin{l1}\label{l-4.2}
If the Goursat problem $\eqref{eq: 2.1}$ and \eqref{eq: 3.10} admits a classical solution in $\Omega_{\tilde{\tau}}$ for some $\tilde{\tau}\in(\tau_0, \infty)$. Then we have 
\begin{align}
    (\bar{\partial}_+\rho, \bar{\partial}_-\rho)\in (M_1(\tau), 0)\times (M_1(\tau), 0), \left(\dfrac{\rho^n\bar{\partial}_+\rho}{\sin^2\delta}, \dfrac{\rho^n\bar{\partial}_-\rho}{\sin^2\delta}\right)\in (M_2(\tau), 0)\times (M_2(\tau), 0)
\end{align}
in $\Omega_{\tilde{\tau}}$.
\end{l1}
\begin{proof}
We first prove that if $0<\delta<\pi/2$ in $\Omega_{\tilde{\tau}}$ then $(\bar{\partial}_+\rho, \bar{\partial}_-\rho)\in (M_1(\tau), 0)\times (M_1(\tau), 0)$ in $\Omega_{\tilde{\tau}}$. We claim that if $(\bar{\partial}_+\rho, \bar{\partial}_-\rho)\in (M_1(\tau), 0)\times (M_1(\tau), 0)$ for all points in $\mathcal{D}_E/\{E\}$ then  $(\bar{\partial}_+\rho, \bar{\partial}_-\rho)\in (M_1(\tau), 0)\times (M_1(\tau), 0)$ at $E$ where $\mathcal{D}_E$ is same as defined in Proposition \ref{p-1}. If $\bar{\partial}_-\rho(E)=M_1(\tau)$ and $\bar{\partial}_+\rho(E)\in [M_1(\tau),0)$. Then by $(\bar{\partial}_+\rho, \bar{\partial}_-\rho)\in (M_1(\tau), 0)\times (M_1(\tau), 0)$ for all $\mathcal{D}_E/\{E\}$ we have $\bar{\partial}_-\bar{\partial}_+\rho\leq 0$ at $E$. However, by second equation of \eqref{eq: 2.15} at the point $E$, we have
\begin{align*}
    c\bar{\partial}_-\bar{\partial}_+\rho&>\dfrac{\tau^4p''(\tau)f}{4c\cos^2\delta}M_1^2(\tau)+M_1(\tau)\sin 2\delta\\
    &>\dfrac{\tau^4p''(\tau)\sin^2 \delta}{2c\cos^2\delta}M_1^2(\tau)+M_1(\tau)\sin 2\delta>0,
\end{align*}
which leads to a contradiction. Hence by a continuity argument we must have $(\bar{\partial}_+\rho, \bar{\partial}_-\rho)\in (M_1(\tau), 0)\times (M_1(\tau), 0)$ in $\Omega_{\tilde{\tau}}$. 

Now for the other part of Lemma, we need to prove that if $0<\delta<\pi/2$ in $\Omega_{\tilde{\tau}}$ then $\left(\dfrac{\rho^n\bar{\partial}_-\rho}{\sin^2\delta}\right)\in (M_2(\tau), 0)\times (M_2(\tau), 0)~~\mathrm{for}~\tau\in[\tau_0,\infty)$ in $\Omega_{\tilde{\tau}}$. Similar to the first part of proof we assume a point $E$ in $\Omega_{\tilde{\tau}}$ such that $\left(\dfrac{\rho^n\bar{\partial}_-\rho}{\sin^2\delta}\right)(E)=M_2(\tau)$ and $\left(\dfrac{\rho^n\bar{\partial}_+\rho}{\sin^2\delta}\right)(E)\in [M_2(\tau),0)$. Then if $\left(\dfrac{\rho^n\bar{\partial}_+\rho}{\sin^2\delta},\dfrac{\rho^n\bar{\partial}_-\rho}{\sin^2\delta}\right)\in(M_2(\tau), 0)\times (M_2(\tau), 0)$ for all points in $\mathcal{D}_E/\{E\}$ then $\bar{\partial}_-\left(\dfrac{\rho^n\bar{\partial}_+\rho}{\sin^2\delta}\right)\leq 0$ at $E$. However noting that $\mathcal{G}>0$ for sufficiently large $n$ we have $\bar{\partial}_-\left(\dfrac{\rho^n\bar{\partial}_+\rho}{\sin^2\delta}\right)> 0$ at $E$, which leads to a contradiction. Thus we have $\left(\dfrac{\rho^n\bar{\partial}_+\rho}{\sin^2\delta}\right)(E)> M_2(\tau)$. Hence we have proved the Lemma.
\end{proof}
It is easy to see that the Lemma \ref{l-4.1} and Proposition \ref{p-1} can be extended to the whole interaction domain $\Omega_{\tilde{\tau}}$ for any $\tilde{\tau}\in [\tau_0, \infty)$ if the curve $\bar{\delta}(\tau)$ possesses countable number of points of extrema in $[\tau_0, \infty)$ and $\underset{\tau\longrightarrow\infty}{\lim}\bar{\delta}(\tau)=\bar{\delta}^*$ exists. Let $\bigcup_{i\in \mathbb{N}}\tau_{i-1}$ is a countable collection of points of local extrema of the curve $\bar{\delta}(\tau)$ in $[\tau_0, \infty)$. Then we can repeat the process of Lemma \ref{l-4.1} and Proposition \ref{p-1} and use Lemma \ref{l-4.2} to obtain the following Lemma.
\begin{l1}
If the Goursat problem \eqref{eq: 2.1} and \eqref{eq: 3.10} admits a $C^1$ solution on $\Omega_{\tilde{\tau}}$ where $\tau\in [\tau_0, \infty)$ and  $\bigcup_{i\in \mathbb{N}}\tau_{i-1}$ is a countable collection of points of local extrema of the curve $\bar{\delta}(\tau)$ in $[\tau_0, \infty)$. Then there exists a sufficiently small $\epsilon_2>0$ which depends only on $\tau$ such that $(\alpha, \beta)$ lies in the invariant region
\begin{align*}
 \Gamma^*:\left(-\dfrac{\pi}{2}-\epsilon_1+2\bar{\delta}(\tau_0)+\chi_\infty, \alpha_0+\epsilon_1+\psi_\infty\right)\times \left(-\dfrac{\pi}{2}-\epsilon_1-\psi_\infty, \alpha_0+\epsilon_1-2\bar{\delta}(\tau_0)-\chi_\infty\right)
    \cap \{\alpha-\beta>\epsilon_2\} ~~\mathrm{in}~\Omega_{\tilde{\tau}}
    \end{align*}
where
\begin{align}
    \psi_\infty=\bar{\delta}^*-\bar{\delta}(\tau_0)+\displaystyle{\int_{\tau_0}^{\infty}}|\bar{\delta}'(\tau)|d\tau=\begin{cases}
   \sum_{i\in \mathbb{N}} 2(\bar{\delta}(\tau_{2i-1})-\bar{\delta}(\tau_{2i-2})), \bar{\delta}(\tau_0)<\bar{\delta}(\tau_1),\\
   \sum_{i\in \mathbb{N}} 2(\bar{\delta}(\tau_{2i})-\bar{\delta}(\tau_{2i-1})), \bar{\delta}(\tau_0)>\bar{\delta}(\tau_1),
    \end{cases},\\
     \chi_\infty=\bar{\delta}^*-\bar{\delta}(\tau_0)-\displaystyle{\int_{\tau_0}^{\infty}}|\bar{\delta}'(\tau)|d\tau=\begin{cases}
   \sum_{i\in \mathbb{N}} 2(\bar{\delta}(\tau_{2i})-\bar{\delta}(\tau_{2i-1})), \bar{\delta}(\tau_0)<\bar{\delta}(\tau_1),\\
      \sum_{i\in \mathbb{N}} 2(\bar{\delta}(\tau_{2i-1})-\bar{\delta}(\tau_{2i-2})), \bar{\delta}(\tau_0)>\bar{\delta}(\tau_1).
    \end{cases}
\end{align}
\end{l1}
\begin{figure}
    \centering
    \includegraphics[width=3 in]{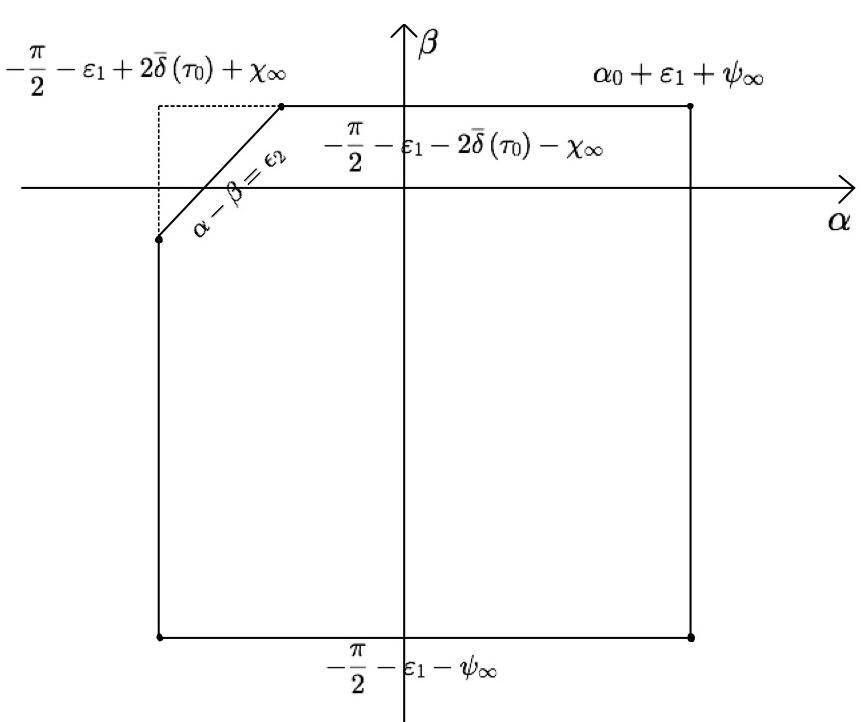}
    \caption{Invariant pentagon}
    \label{fig: 9}
\end{figure}
\begin{proof}
By repeating the steps of Proposition \ref{p-1} for $\tau=\tau_0, \tau_1, \tau_2,.......,\tau_n,\tau_{n+1},......$ it is easy to see that for all $\tau\in[\tau_0, \infty)$ we have
$$\left(\alpha, \beta\right)\in \left(-\dfrac{\pi}{2}-\epsilon_1+2\bar{\delta}(\tau_0)+\chi_\infty, \alpha_0+\epsilon_1+\psi_\infty\right)\times \left(-\dfrac{\pi}{2}-\epsilon_1-\psi_\infty, \alpha_0+\epsilon_1-2\bar{\delta}(\tau_0)-\chi_\infty\right)$$
in $\Omega_{\tilde{\tau}}.$

Now we move forward to prove the other part of the Lemma. Now if $0<\alpha-\beta=\epsilon_2<\alpha_0+\dfrac{\pi}{2}$ at some point in $\Omega_{\tilde{\tau}}$ then we see from \eqref{eq: 2.12} that
\begin{align*}
    c\bar{\partial}_\pm \left(\dfrac{\alpha-\beta}{2}\right)&=\dfrac{\sin^2\delta}{2}\bigg[2+\dfrac{p''(\tau)}{4c\rho^{4+n}}\bigg[\tan\delta-\dfrac{\Omega\sin2\delta}{2}\bigg]\left(\dfrac{\rho^n\bar{\partial}_\pm\rho}{\sin^2\delta}\right)\bigg],\\
    &>\dfrac{\sin^2\delta}{2}\bigg[2+\dfrac{p''(\tau)}{4c\rho^{4+n}}\bigg[\tan\delta-\dfrac{\Omega\sin2\delta}{2}\bigg]M_2(\tau)\bigg]>0.
\end{align*}
Since $(\alpha-\beta)(P)=\alpha_0+\pi/2>\epsilon_2$ then we have $\alpha-\beta>\epsilon_2$ in $\Omega_{\tilde{\tau}}$. Hence proved.
\end{proof}
The above Lemma can be generalized for a more general case of $\bar{\delta}(\tau)$ as follows:
\begin{l1}\label{l-4.4}
\textit{Assume that the Goursat problem $\eqref{eq: 2.1}$, $\eqref{eq: 3.10}$ admits a $C^1$ solution in the domain $\Omega_{\Tilde{\tau}}$ then the invariant region for $(\alpha, \beta)$ is
\begin{align*}
 \left(-\dfrac{\pi}{2}-\epsilon_1+2\bar{\delta}(\tau_0)+\chi({\tilde{\tau}}), \alpha_0+\epsilon_1+\psi({\tilde{\tau}})\right)\times \left(-\dfrac{\pi}{2}-\epsilon_1-\psi({\tilde{\tau}}), \alpha_0+\epsilon_1-2\bar{\delta}(\tau_0)-\chi({\tilde{\tau}})\right)
 \cap\{\alpha-\beta>\epsilon_2\}. 
    \end{align*}}
\end{l1}

From the above Lemmas, it follows that if the Goursat problem \eqref{eq: 2.1} and \eqref{eq: 3.10} admits a $C^1$ solution on $\Omega_{\tilde{\tau}}$ for any $\tilde{\tau}\in [\tau_0, \infty)$ then we have
\begin{align}\label{eq: 74}
    0<\dfrac{\epsilon_2}{2}<\delta<\epsilon_1+\psi_{\infty}<\dfrac{\pi}{2}~~\mathrm{on}~\Omega_{\Tilde{\tau}}.
\end{align}
\subsection{$C^0$ and $C^1$ norm estimates of the solution}
\begin{l1}{($C^0$ estimates)}\label{l-4.5}
If the Goursat problem \eqref{eq: 2.1} and \eqref{eq: 3.10} admits a $C^1$ solution in $\Omega_{\tilde{\tau}}$, where $\tilde{\tau}\in [\tau_0, \infty)$. Then there exists a positive $\mathcal{H}_0$ depending on $\tilde{\tau}$ such that $||(u, v, \rho)||_{C^0(\Omega_{\Tilde{\tau}})}<\mathcal{H}_0.$
\end{l1}
\begin{proof}
In view of the relations \eqref{eq: 2.7} and Lemma \ref{l-4.4}, this Lemma can be easily proved, so we omit the details.
\end{proof}
\begin{l1}{($C^1$ estimates)}\label{l-4.6}
Suppose that the Goursat problem \eqref{eq: 2.1} and \eqref{eq: 3.10} admits a $C^1$ solution in $\Omega_{\Tilde{\tau}}$ for some $\tilde{\tau}\in(\tau_0, \infty)$. Then there exists a positive $\mathcal{H}_1$ depending on $\tilde{\tau}$ such that $||(u, v, \rho)||_{C^1(\Omega_{\Tilde{\tau}})}<\mathcal{H}_1.$
\end{l1}
\begin{proof}
This result can be immediately obtained using \eqref{eq: 2.13}, \eqref{eq: 74}, Lemma \ref{l-4.2}, \ref{l-4.4}, \ref{l-4.5} and the relations
$$\partial_\xi=\dfrac{\sin\alpha\bar{\partial}_--\sin\beta\bar{\partial}_+}{\sin 2\delta},~~\partial_\eta=\dfrac{\cos\beta\bar{\partial}_+-\cos\alpha\bar{\partial}_-}{\sin 2\delta}.$$
\end{proof}
\section{Existence of global solution}\label{5}
First let us assume that the Goursat problem $\eqref{eq: 2.1}$, $\eqref{eq: 3.10}$ admits a unique $C^1$ solution in $\Omega_{\tilde{\tau}}$ where $\tilde{\tau}\in [\tau_0, \infty)$. Then by using Lemma \ref{l-4.4}, one can easily see that along the level curves of specific volume $\tau$
\begin{align}\label{eq: 59}
    \bigg|\dfrac{d\xi}{d\eta}\bigg|&=\bigg|\dfrac{\cos\beta\bar{\partial}_+\rho-\cos\alpha\bar{\partial}_-\rho}{\sin\beta\bar{\partial}_+\rho-\sin\alpha\bar{\partial}_-\rho}\bigg|\nonumber\\
    &<\bigg|\dfrac{\cos\beta\bar{\partial}_+\rho}{\sin\beta\bar{\partial}_+\rho-\sin\alpha\bar{\partial}_-\rho}\bigg|+\bigg|\dfrac{\cos\alpha\bar{\partial}_-\rho}{\sin\beta\bar{\partial}_+\rho-\sin\alpha\bar{\partial}_-\rho}\bigg|\nonumber\\
    &<\bigg|\dfrac{1}{\sin\beta}\bigg|+\bigg|\dfrac{1}{\sin\alpha}\bigg|<\dfrac{2}{\sin(2\bar{\delta}(\tau_0)+\chi({\tilde{\tau}})-\epsilon_1)}.
\end{align}
Then we are going to prove the existence of a global solution by extending the local solution in this Section. For this purpose, we solve several local Goursat problems in each extension step. Let $\overline{XZ}$ and $\overline{XY}$ be positive($C+$) and negative($C-$) characteristics in $\Omega_{\tilde{\tau}}$, respectively, where $Y=(\xi_Y, \eta_Y)$ and $Z=(\xi_Z, \eta_Z)$ are two points lying on the level curve $\tau=\tilde{\tau}$. Then, we prescribe the boundary data
\begin{align}\label{eq: 5.1}
    (u, v, \tau)=\begin{cases}
    (u|_{XZ}, v|_{XZ}, \tau|_{XZ}), ~~~~~~~~~~~on~~~~~~ \overline{XZ}\\
    (u|_{XY}, v|_{XY}, \tau|_{XY}), ~~~~~~~~~~~on~~~~~~ \overline{XY}
    \end{cases}
\end{align}
where ($u|_{XY}, v|_{XY}, \tau|_{XY}$) and ($u|_{XZ}, v|_{XZ}, \tau|_{XZ}$)  are the values of $(u, v, \tau)$ on $\overline{XY}$ and $\overline{XZ}$, respectively.

Then we have the following Lemma.
\begin{l1}\label{l-5.1}
\textit{ If $\eta_{Z}-\eta_{Y}<-\dfrac{\sin(2\bar{\delta}(\tau_0)+\chi(\tilde{\tau})-\epsilon_1)}{2\tilde{\tau}M_1(\tilde{\tau})}$,
then the Goursat problem $\eqref{eq: 2.1}$, $\eqref{eq: 5.1}$ admits a global $C^1$ solution in a domain bounded by $\overline{XY}, \overline{XZ}, \overline{TY}$ and $\overline{TZ}$ where $\overline{TY}$ is the $C_+$ characteristic passing through $Y$ while $\overline{TZ}$ is the $C_-$ characteristic passing through $Z$. Moreover, this solution satisfies $\tilde{\tau}>\dfrac{\tau}{2}$.}  
\end{l1}
\begin{proof} 
Let $K$ be any point in the quadrilateral domain bounded by $\overline{XY}, \overline{XZ}, \overline{TY}$ and $\overline{TZ}$ such that the $C_+$ characteristics passing through $K$ intersects with $XY$ at $K_+$ while $C_-$ characteristics passing through $K$ intersects with $XZ$ at $K_-$.Then using Lemma \ref{l-4.2} we have along $KK_-$
\begin{align*}
    \dfrac{1}{\tau(K)}&=\dfrac{1}{\tau(K_-)}+\displaystyle{\int_{K_-K}\bar{\partial}_-\rho}\\
    &>\dfrac{1}{\tau(K_-)}+\dfrac{M_1(\tau)(\eta_Y-\eta_Z)}{\sin(2\bar{\delta}(\tau_0)+\chi(\tilde{\tau})-\epsilon_1)}>\dfrac{1}{2\tilde{\tau}}.
\end{align*}
Hence, using \eqref{eq: 2.13} we obtain a $C^1$ norm estimate of the solution to the Goursat problem \eqref{eq: 2.1} and \eqref{eq: 3.10}, so by the theory of global classical solution for quasilinear hyperbolic equations \cite{li1985boundary}, the Lemma can be proved.
\end{proof}
\begin{figure}
    \centering
    \includegraphics[width= 4.7 in]{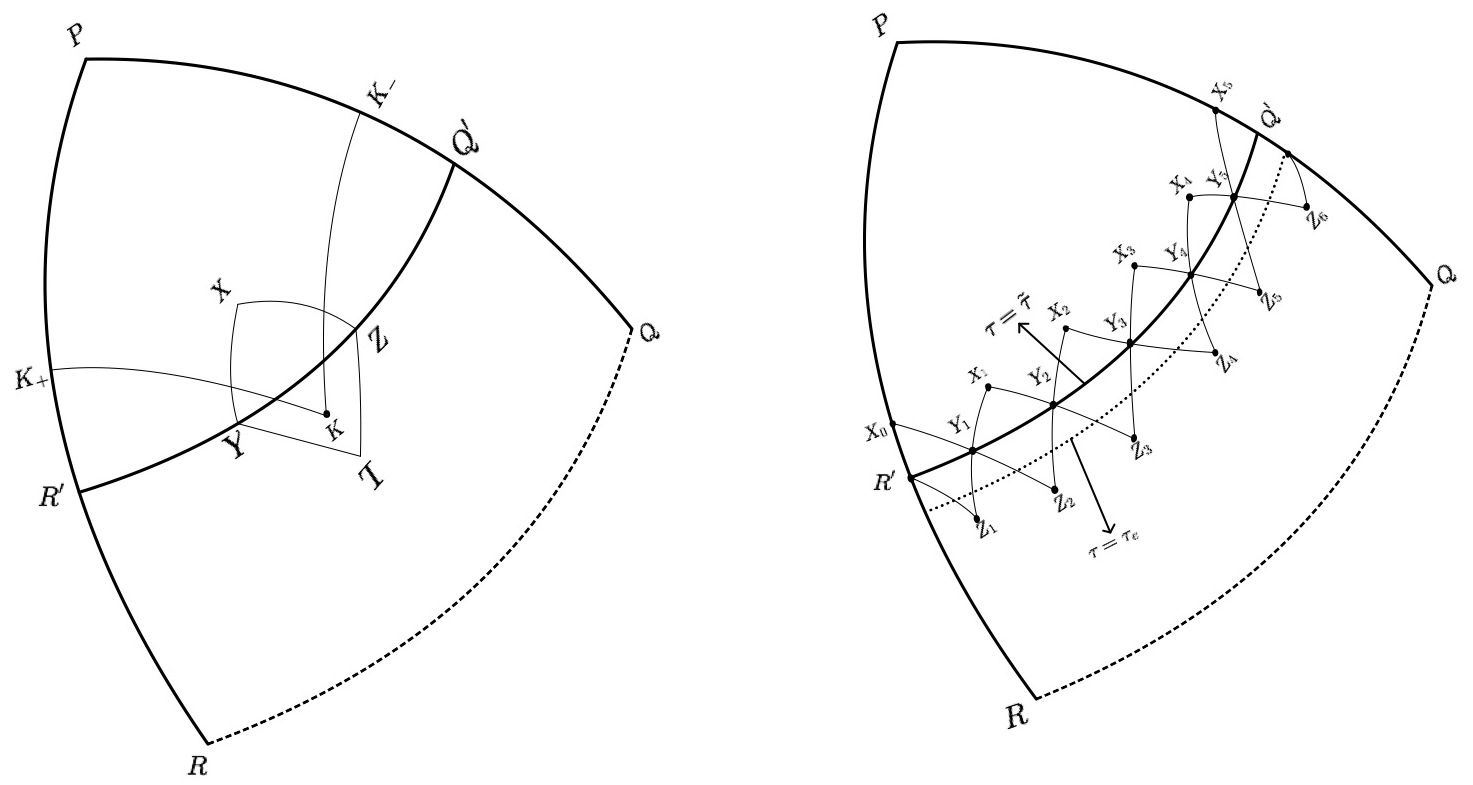}
    \caption{Extension of local solution to global solution}
    \label{fig:my_label}
\end{figure}
The main result of this work can be summarized in the following theorem:

\begin{t1}\label{t-5.1}
\textit{(Global solution and regularity of vacuum boundary)~If the wall angle satisfies $\theta\in (-\pi/2, 0)$ and $2\bar{\delta}(\tau_0)+\chi(\tau)<\alpha_0+\dfrac{\pi}{2}<4\bar{\delta}(\tau_0)$ then the Goursat problem \eqref{eq: 2.1} with the boundary data $\eqref{eq: 3.10}$ admits a unique global $C^1$ solution in the region $PQR$ where the curve $\overline{QR}$ is an interface between gas and vacuum connecting $Q$ and $R$, which is Lipschitz continuous.}
%\begin{figure}
%\centering
 %   \includegraphics[width= 7 in]{Extension3.eps}
  %  \caption{Left: A curved quadrilateral building block; Right: Global %solution}
   % \label{fig:my_label5}
%\end{figure}
\end{t1}
\begin{proof}
Let $Y_0=Q', Y_1, Y_2, Y_3,........, Y_n=R'$ be $n+1$ different points on the level curve $\tau=\tilde{\tau}$ such that
$$0<\eta_{Y_{i+1}}-\eta_{Y_{i}}<-\dfrac{\sin(2\bar{\delta}(\tau_0)+\chi(\tilde{\tau})-\epsilon_1)}{2\tilde{\tau}M(\tilde{\tau})}$$
Then we draw a $C_+$ characteristic curve from the point $Y_i$ which cuts the $C_-$ characteristic curve passing through $Y_{i+1}$ at a point $X_i$ for $i=0,1, 2,....., n-1$. Therefore, noting the fact that $\bar{\partial}_\pm \tau > 0$ it is easy to see that the level curve $\tau=\tilde{\tau}$ is a non-characteristic curve which means that $X_i\neq Y_i$ and $X_i\neq Y_{i+1}$ for any $i=0, 1, 2,......, n-1$.

From the point $X_i$, we build a small Goursat problem below the level curve $\tau=\tilde{\tau}$ where $\overline{X_iY_i}$ and $\overline{X_i Y_{i+1}}$ are the characteristic boundaries in the domain $\Omega_{\tilde{\tau}}$. Therefore, we can use Lemma $\ref{l-5.1}$ to conclude that the Goursat problem \eqref{eq: 2.1}, \eqref{eq: 5.1} with the characteristic boundaries $\overline{X_i Y_i}$ and $\overline{X_i Y_{i+1}}$ admits a $C^1$ solution in the quadrilateral domain bounded by $\overline{X_i Y_i}, \overline{X_i Y_{i+1}}, \overline{Y_i Z_{i+1}}$ and $\overline{Y_{i+1} Z_{i+1}}$ where $\overline{Y_i Z_{i+1}}$ is the $C_-$ characteristic curve passing through $Y_i$ and $\overline{Y_{i+1} Z_{i+1}}$ is the $C_+$ characteristic curve passing through $Y_{i+1}$.

Assume that $Z_0=Q$ and $Z_{n+1}=R$, then for each $i=0, 1, 2,.....,n$, there exists a $\tau_i\in(\tau_0, \infty)$, such that the Goursat problem for $\eqref{eq: 2.1}$ admits a unique $C^1$ solution in the domain closed by $\overline{Y_i Z_{i}}$, $\overline{Y_i Z_{i+1}}$ and the level curve $\tau(\xi, \eta)=\tau_i$ with $\overline{Y_i Z_i}$ and $\overline{Y_i Z_{i+1}}$ as the characteristic boundaries. If we denote $\tau_e=\min\{\tau_0, \tau_1, \tau_2,......., \tau_n, \tau(Z_1), \tau(Z_2),......,\tau(Z_n)\}$ then utilizing the fact $\bar{\partial}_+ \tau>0$ we see that $\tau_e>\tilde{\tau}$. Then, the solution is extended from $\Omega_{\tilde{\tau}}$ to $\Omega_{\tau_e}$. Repeating the above process, one can construct the global solution of the Goursat problem in the whole domain $PQR$. 

Further, from \eqref{eq: 59} we see that the family of level curves of specific volume are Lipschitz continuous, hence we can use Arzela-Ascoli's theorem to conclude that the vacuum boundary $\overline{QR}$ is also Lipschitz continuous.
\end{proof}

\section{Applications}\label{6}
In this section, we first formulate two-dimensional modified shallow water equations which can be considered as a particular case of two-dimensional Euler equations with a special equation of state. Then we obtain a global solution to a dam-break type problem for the two-dimensional modified shallow water equations as one of the applications of our work. Further, we consider some special equations of states and solve the gas expansion problem through a sharp corner for them while recovering some of the results from the available literature.
\subsection{Two-constant equation of state}
Let us consider an equation of state of the form \cite{hu2014axisymmetric}
\begin{align}
    p(\tau)=A_1\tau^{\gamma_1}+B_1\tau^{\gamma_2}
\end{align}
where $\gamma_1, \gamma_2$, $A_1$ and $B_1$ are constants. Such kind of equations of states are relevant in many physical models, for example, it can be taken as the sum of the fluid and magnetic pressure in the ideal magnetogasdynamics where the magnetic field is orthogonal to the velocity vector \cite{sekhar2010riemann}. Also, it can be used to express two-dimensional modified shallow water equations \cite{sekhar2008interaction} where the gas expansion problem can be expressed as a dam-break type problem. One other example of such an equation of state is extended Chaplygin gas which has been utilized in the recent past by many mathematicians and physicists as a candidate to explain the accelerated expansion of the universe \cite{pourhassan2014extended}. We consider the first two cases one by one in this subsection and solve the corresponding 2-D Riemann problem for them.
\subsubsection{Modified shallow water equations in two-dimensions}
As an application of the gas expansion problem, we consider a dam-break type problem for two-dimensional modified shallow water equations and construct a global solution for the dam-break type problem using the construction technique of this article. We first formulate two-dimensional modified shallow water equations using the ideas from \cite{karelsky2006particular}. We consider isentropic incompressible Euler equations in three dimensions for a free surface fluid of the form
\begin{align}\label{eq: 64}
\begin{cases}
    \nabla.\vec{\mathbf{u}}=0,\\
    \vec{\mathbf{u}}_t+\vec{\mathbf{u}}\cdot\nabla.\vec{\mathbf{u}}+\dfrac{1}{\rho}\nabla p+\vec{\mathbf{g}}=0,
    \end{cases}
\end{align}
where $\vec{\mathbf{u}}=(u, v, w)$ is the velocity of the fluid, $p$ is the pressure and $\vec{\mathbf{g}}=(0,0, g)$ is the gravitational acceleration. We apply the no penetration boundary condition at the lower boundary while the kinematic boundary condition at the free surface so that
\begin{align}\label{eq: 63}
    \begin{cases}
        w|_{z=0}=0,\\
        w|_{z=h}=\dfrac{\partial h}{\partial t}+u\dfrac{\partial h}{\partial x}+v\dfrac{\partial h}{\partial y}.
    \end{cases}
\end{align}
Further, we assume a hydrostatic pressure field of the form $p=p_0+gh(h-z)$, where $h(x, y)$ is the fluid depth and $p_0$ is the atmospheric pressure. Therefore, $z$-momentum equation from \eqref{eq: 64} implies that 
\begin{align}
    \dfrac{Dw}{Dt}=\dfrac{\partial w}{\partial t}+u\dfrac{\partial w}{\partial x}+v\dfrac{\partial w}{\partial y}+w\dfrac{\partial w}{\partial z}=0,
\end{align}
which reduces the governing system as
\begin{align}\label{eq: 66}
\begin{cases}
    \dfrac{\partial u}{\partial x}+ \dfrac{\partial v}{\partial y}+ \dfrac{\partial w}{\partial z}=0,\\
    \dfrac{\partial u}{\partial t}+ \dfrac{\partial u^2}{\partial x}+ \dfrac{\partial uv}{\partial y}+\dfrac{\partial uw}{\partial z}+\dfrac{1}{\rho}\dfrac{\partial p}{\partial x}=0,\\
    \dfrac{\partial v}{\partial t}+ \dfrac{\partial uv}{\partial x}+ \dfrac{\partial v^2}{\partial y}+\dfrac{\partial vw}{\partial z}+\dfrac{1}{\rho}\dfrac{\partial p}{\partial y}=0.
    \end{cases}
\end{align}
Finally we can integrate equations in \eqref{eq: 66} with the boundary conditions \eqref{eq: 63} to obtain 
\begin{align}\label{eq: 67}
    \begin{cases}
       \dfrac{\partial h}{\partial t}+\dfrac{\partial }{\partial x}\displaystyle{\int_0^{h}udz}+\dfrac{\partial }{\partial y}\displaystyle{\int_0^{h}vdz}=0,\\
       \dfrac{\partial }{\partial t}\displaystyle{\int_0^{h}udz}+\dfrac{\partial }{\partial x}\displaystyle{\int_0^{h}u^2dz}+\dfrac{\partial }{\partial y}\displaystyle{\int_0^{h}uvdz}+gh\dfrac{\partial h}{\partial x}=0,\\
       \dfrac{\partial }{\partial t}\displaystyle{\int_0^{h}vdz}+\dfrac{\partial }{\partial x}\displaystyle{\int_0^{h}uvdz}+\dfrac{\partial }{\partial y}\displaystyle{\int_0^{h}v^2dz}+gh\dfrac{\partial h}{\partial y}=0.
    \end{cases}
\end{align}
Let us represent the horizontal velocities $u$ and $v$ as a sum of depth-averaged velocity and its variation from mean velocity such that $u=\bar{u}+u^*, v=\bar{v}+v^*$, where $\bar{u}=\dfrac{1}{h}\displaystyle{\int_{0}^{h}udz}$ and $\displaystyle{\int_{0}^{h}u^*dz}=0$ while $\bar{v}=\dfrac{1}{h}\displaystyle{\int_{0}^{h}vdz}$ and $\displaystyle{\int_{0}^{h}v^*dz}=0$.

Then we can simplify the system \eqref{eq: 67} to obtain
\begin{align}
    \begin{cases}
        h_t+(h \bar{u})_x+(h \bar{v})_y=0,\\
        (h \bar{u})_t+(h \bar{u}^2+\dfrac{gh^2}{2})_x+(h \bar{u}\bar{v})_y+R^1_x+R^2_y=0,\\
        (h \bar{v})_t+(h \bar{u}\bar{v})_x+(h \bar{v}^2+\dfrac{gh^2}{2})_y+R'^1_y+R'^2_x=0
    \end{cases}
\end{align}
where $R^1=\displaystyle{\int_0^{h}u^{*2}dz}$, $R^2=\displaystyle{\int_0^{h}u^{*}v^*dz}$, $R'^1=\displaystyle{\int_0^{h}v^{*2}dz}$ and $R'^2=\displaystyle{\int_0^{h}u^* v^{*}dz}$ are the nonlinear terms describing the effect of advective transport of the impulse caused by velocity difference. For simplicity, we consider the case where $u^{*2}\approx v^{*2}\approx k_0$ as in \cite{karelsky2006particular} and assume that the effect of the nonlinear term $u^*v^*$ is negligible so that the modified shallow water equations in two-dimensions can be obtained of the following form
\begin{align}\label{eq: 72}
    \begin{cases}
        h_t+(h \bar{u})_x+(h \bar{v})_y=0,\\
        (h \bar{u})_t+(h \bar{u}^2+\dfrac{gh^2}{2}+kh)_x+(h \bar{u}\bar{v})_y=0,\\
        (h \bar{v})_t+(h \bar{u}\bar{v})_x+(h \bar{v}^2+\dfrac{gh^2}{2}+kh)_y=0,
    \end{cases}
\end{align}
where $k=k_0/g$ is a reduced factor characterizing advective transport of the impulse. 

It is worth noting that the above modified shallow water equations are comparable with Euler equations in two dimensions with an equation of state of the form
$p(\tau)=\dfrac{A_1}{\tau}+\dfrac{B_1}{\tau^2}$ where $A_1=k$, $B_1=g/2$ and $\tau=1/h$. Hence we directly calculate
\begin{align*}
    &p(\tau)=A_1\tau^{-1}+B_1\tau^{-2},~~~p'(\tau)=-A_1\tau^{-2}-2B_1\tau^{-3}<0,~~~p''(\tau)=2A_1\tau^{-3}+6B_1\tau^{-4}>0,\\
    &m(\tau)=\dfrac{A_1\tau^{-2}+B_1\tau^{-3}}{A_1\tau^{-2}+3B_1\tau^{-3}}>0,~~~ m'(\tau)=\dfrac{2B_1A_1\tau^{-6}}{[A_1\tau^{-2}+3B_1\tau^{-3}]^2}>0.
\end{align*}
Hence, we can use the Lemma \ref{l-4.1}, Proposition \ref{p-1} and Theorem \ref{t-5.1} to obtain the following result for modified shallow water equations in two dimensions.
\begin{t1}
If the wall angle satisfies $\theta\in (-\pi/2, 0)$ and $2\bar{\delta}(\tau_0)<\alpha_0+\pi/2<4\bar{\delta}(\tau_0)$ then the modified shallow water system \eqref{eq: 72} with initial data \eqref{eq: 2} admits a global classical $C^1$ solution.
\end{t1}

\subsubsection{Isentropic magnetogasdynamics system in two-dimensions}
Let us consider the following system of isentropic magnetogasdynamics \cite{cabannes1970theoretical}
\begin{align}\label{eq: 83}
    \begin{cases}
        \rho_t+\nabla.(\rho\vec{\mathbf{u}})=0,\\
        (\rho\vec{\mathbf{u}})_t+\nabla.(\rho \vec{\mathbf{u}}\otimes \vec{\mathbf{u}}+p\mathbf{I})-\mu(\nabla\times \vec{\mathbf{H}})\times\vec{\mathbf{H}}=0,\\
        \vec{\mathbf{H}}_t-\nabla\times(\vec{\mathbf{u}}\times\vec{\mathbf{H}})=0,\\
        \nabla.\vec{\mathbf{H}}=0,
    \end{cases}
\end{align}
where $\rho$ is the fluid density, $p=A_1\rho^{\gamma}$ is the fluid pressure($A_1>0$), $\vec{\mathbf{u}}$ is the fluid velocity and $\vec{\mathbf{H}}$ is magnetic field vector.

If we assume that the velocity vector and magnetic field vector are perpendicular to each other or in other words if $\vec{\mathbf{u}}=(u, v, 0)$ and $\vec{\mathbf{H}}=(0, 0, H)$ then we can reduce the system \eqref{eq: 83} as follows (see \cite{chen2020rarefaction} for a complete derivation)
\begin{align}\label{eq: 84}
    \begin{cases}
        \rho_t+(\rho u)_x+(\rho v)_y=0,\\
(\rho u)_t+(\rho u^2+p+\dfrac{\mu}{2}H^2)_x+(\rho u v)_y=0,\\
(\rho v)_t+(\rho u v)_x+(\rho v^2+p+\dfrac{\mu}{2}H^2)_y=0,\\
H_t+(H u)_x+(H v)_y=0.
    \end{cases}
\end{align}
From the first and last equation of system \eqref{eq: 84} it is easy to see that $H=\kappa_0\rho$ along streamlines of the flow which is usually referred to as frozen law in the literature; see viz. \cite{sekhar2010riemann}. Using this substitution in the governing system \eqref{eq: 84} we can obtain a modified pressure of the form $p(\tau)=A_1\tau^{-\gamma}+B_1\tau^{-2}$, where $B_1=\dfrac{\mu \kappa_0^2}{2}$. Then we have
\begin{align*}
    &p(\tau)=A_1\tau^{-\gamma}+B_1\tau^{-2},~~ p'(\tau)=-A_1\gamma\tau^{-\gamma-1}-2B_1\tau^{-3}<0,~~
    p''(\tau)=A_1\gamma(\gamma+1)\tau^{-\gamma-2}+6B_1\tau^{-4}>0,\\
    &m(\tau)=\dfrac{A_1\gamma(3-\gamma)\tau^{-\gamma-1}+2B_1\tau^{-3}}{A_1\gamma(\gamma+1)\tau^{-\gamma-1}+6B_1\tau^{-3}}>0,~~~ m'(\tau)=\dfrac{8B_1A_1\gamma(\gamma-2)^2\tau^{-\gamma-5}}{[A_1\gamma(\gamma+1)\tau^{-\gamma-1}+6B_1\tau^{-3}]^2}>0.
\end{align*}
For this system, we again use Lemma \ref{l-4.1}, Proposition \ref{p-1} and Theorem \ref{t-5.1} to recover the following result from \cite{chen2020expansion}.
\begin{t1}
If the wall angle satisfies $\theta\in (-\pi/2, 0)$ and $2\bar{\delta}(\tau_0)<\alpha_0+\pi/2<4\bar{\delta}(\tau_0)$ then the initial boundary value problem \eqref{eq: 84} and \eqref{eq: 2} admits a global classical $C^1$ solution.
\end{t1}

\subsection{Van der Waals gas}
Let us consider a van der Waals type equation of state of the form \cite{lai2019global}
\begin{align}
    p(\tau)=\dfrac{S_1}{(\tau-1)^{\gamma+1}}-\dfrac{1}{\tau^2},
\end{align}
where $S_1$ is a positive constant in correspondence with the gas entropy and $\gamma$ is the gas constant lying between $0$ and $1$. It is easy to see that for $\tau>\tau_0$, we have $p'(\tau)<0$ and $p''(\tau)>0$ for sufficiently large $\tau>4$ and for any given $S_1\in (1/4, 81/256)$, when $\gamma$ is sufficiently close to $0$ we have that $m(\tau)>0$ and $m'(\tau)<0$; see \cite{lai2015centered} for more details. Therefore, we can use the second part of the Lemma \ref{l-4.1} and Proposition \ref{p-1} and Theorem \ref{t-5.1} to obtain the following result.
\begin{t1}
If the wall angle satisfies $\theta\in (-\pi/2, 0)$ and $2\bar{\delta}(\tau_0)+\chi(\tau)<\alpha_0+\pi/2<4\bar{\delta}(\tau_0)$ then the initial boundary value problem \eqref{eq: 1} and \eqref{eq: 2} admits a global classical $C^1$ solution for van der Waals gas.
\end{t1}
\section{Conclusions and future scope}\label{7}
Here we established the existence of a global solution to a 2-D Riemann problem for compressible Euler equations with a general convex equation of state. It is a well-known fact that the Euler equations in the self-similar plane is a mixed type system so to maintain the hyperbolicity of Euler equations, we constructed the invariant regions of the characteristic variables using the characteristic decomposition method and hence obtained a global solution to a gas expansion problem by extending the local solution into the whole interaction domain. Further, we formulated two-dimensional modified shallow water equations and obtained a global solution for a dam-break type problem as one of the applications of this work. Also, we recovered some of the existence results of the gas expansion problems for certain equations of states from the available literature. It was worth noting that we considered the convex equation of state and a special assumption on wall angle in this work but in future, we would like to relax the restriction on the equation of state and solve this problem for any arbitrary equation of state has more than one inflection point. Moreover, we also wish to discuss different wall angles which may include interactions of composite waves and shocks also.\\
\textbf{Acknowledgements}~ \textit{The first author (RB) gratefully acknowledges the research support from the University Grant Commission, Government of India.  The second author (TRS) would like to
thank SERB, DST, Government of India (Ref. No. CRG/2022/006297) for its financial support through
Core research grant.}
\biboptions{sort&compress}
\bibliographystyle{elsarticle-num}
\bibliography{Citation}
\end{document}